\documentclass[mnsc,nonblindrev]{informs3arXiv} 
\pdfoutput=1	
\usepackage{amsmath}
\usepackage{amssymb}
\usepackage{multirow}
\usepackage{setspace}
\usepackage{algorithm}
\usepackage{algorithmic}
\usepackage{pdflscape}
\usepackage{setspace}
\usepackage{subfigure}
\usepackage{authblk}
\usepackage{amsfonts}	
\usepackage{xcolor}
\usepackage{color}
\usepackage{multirow}
\usepackage{graphicx}
\usepackage[small,compact]{titlesec}
\usepackage{calc}
\usepackage{verbatim}
\usepackage{mathtools}
\usepackage{tikz}
\usepackage{float}
\usepackage{booktabs}
\usetikzlibrary{shapes.geometric}

\newcommand{\bh}{ \mathbf{h} }

\newcommand{\bX}{ \mathbf{X} }
\newcommand{\bx}{ \mathbf{x} }
\newcommand{\bY}{ \mathbf{Y} }
\newcommand{\by}{ \mathbf{y} }

\newcommand{\bz}{ \mathbf{z} }

\OneAndAHalfSpacedXI 

\usepackage{natbib}
 \bibpunct[, ]{(}{)}{,}{a}{}{,}%
 \def\bibfont{\small}%
\TheoremsNumberedThrough     

\EquationsNumberedThrough    

\begin{document}


\RUNAUTHOR{Boutilier and Chan}

\RUNTITLE{Drone-delivered AEDs}

\TITLE{Response time optimization for drone-delivered automated external defibrillators}

\ARTICLEAUTHORS{%
\AUTHOR{Justin J. Boutilier, Timothy C.Y. Chan}
\AFF{Department of Mechanical and Industrial Engineering, University of Toronto, Toronto, Ontario M5S 3G8, Canada \\ \EMAIL{j.boutilier@mail.utoronto.ca} \EMAIL{ tcychan@mie.utoronto.ca}}
} 

\ABSTRACT{%
Out-of-hospital cardiac arrest (OHCA) claims over $400,000$ lives each year in North America and is one of the most time-sensitive medical emergencies. Drone-delivered automated external defibrillators (AEDs) have the potential to be a transformative innovation in the provision of emergency care for OHCA. In this paper, we propose a simulation-optimization framework to minimize the total number of drones required to meet a pre-specified response time goal, while guaranteeing a sufficient number of drones are located at each base. To do this, we develop a location-queuing model that is based on the $p$-median architecture, where each base constitutes an explicit $M/M/d$ queue, and that incorporates estimated baseline response times to the demand points. We then develop a reformulation technique that exploits the baseline response times, allowing us to solve real-world instances to optimality using an off-the-shelf solver. To test our model, we develop a two-stage machine learning approach to simulate both the locations and baseline response times for future OHCAs. We demonstrate the application of our framework using eight years of real data from an area covering 26,000 square kilometers around Toronto, Canada.  A modest number of drones are required to significantly reduce response times in all regions. Furthermore, an objective function focused on improving the 90th percentile is well-suited for use in practice because the model reduces the entire response time distribution, while providing equitable coverage in both cities and rural areas. Overall, this paper provides a realistic framework that can be leveraged by healthcare providers seeking to implement a drone network.
}%


\KEYWORDS{Facility location, queuing, emergency medical services, drones, optimization, machine learning}
\HISTORY{}

\maketitle

%

\section{Introduction}

Out-of-hospital cardiac arrest (OHCA) claims over $400,000$ lives each year in North America \citep{Go2013}. It is one of the most time-sensitive medical emergencies with survival estimated to decrease by 7-10\% for each minute without treatment \citep{Valenzuela1997}. An automated external defibrillator (AED) is one of the only effective methods for treating OHCA and studies have shown that the probability of survival increases by 40-70\% with prompt defibrillation \citep{Caffrey2002,PAD2004}. As a result, many municipalities have developed public access defibrillation (PAD) programs that place defibrillators in public areas so that they can be used by bystanders during an OHCA event \citep{PAD2004}. Indeed, PAD programs have demonstrated great efficacy for OHCAs in public locations by reducing the time to defibrillation and improving survival outcomes \citep{Hazinski2005}.

While progress has been made towards improving survival for public location OHCAs, the vast majority of OHCAs occur in private, residential settings with slower emergency response times and lower survival, especially in rural areas   . There is a fundamental coverage limit of public AEDs that limits their usefulness for private OHCAs \citep{Siddiq2013}. Moreover, statically deployed public AEDs have low historical utilization \citep{Weisfeldt2010}, which can be attributed in part to access and availability issues \citep{Sun2016}. Thus, a new approach to improving AED access and reducing the time to defibrillation is needed, especially for residential and rural areas.

In the past few years, several organizations and researchers have developed drone technology that can be used to deliver AEDs \citep{Delft, Reece}. Drone-delivered AEDs have received increasing attention due to their potential to be a transformative innovation in the provision of emergency care. However, previous research on drone-delivered AEDs is limited to three studies, all of which are based on covering-type location models. Two of these studies focus solely on optimizing drone base locations without considering system congestion \citep{Pulver2016, Claesson2016}. The third paper considers both optimal base locations and system congestion, where a queuing model to determine the number of drones required at each base is solved after the base locations have been optimized \citep{Boutilier2017}.

While the aforementioned approaches are novel in their application of optimization ideas to drone deployment in a healthcare application, they have several limitations. First, they all use coverage-type models. Modeling for drone-delivered AEDs should ideally be based on a median-style model to optimize response time explicitly. Drones for emergency response will need to be integrated with an emergency medical services (EMS) system and the performance of those systems is typically measured via average and 90th percentile response time metrics. One of the major challenges of dealing with response time models is the computational difficulty in solving large scale instances due to the lack of sparsity in the response time matrix \citep{Li2011}. That is why coverage models are often preferred; sparsity is induced by removing demand points outside the coverage radius of each facility. Second, any consideration of congestion should ideally be combined with the location decision in an integrated model; separate models solved sequentially may generate suboptimal solutions. We expect AED-enabled drones to have high utilization because they would be dispatched to all suspected OHCAs, which underlies the importance of integrating queuing constraints with the location problem. Third, none of the previous models accounted for the baseline response times (i.e., historical 911 response times) in their respective study regions when locating drone resources. Given that this application of drone technology involves layering drones on top of existing emergency response resources, the goal should be to optimize overall system performance. For example, drones could be used to patch geographic holes that are hard to reach or that suffer from slower ambulance response times.


In this paper, we address the aforementioned limitations of existing models in a single optimization framework. Our specific contributions are:

\begin{enumerate}

\item We propose a novel integrated location-queuing model that is based on the $p$-median architecture, where each base constitutes an explicit $M/M/d$ queue, and that incorporates estimated baseline response times to the demand points. We use our framework to minimize the total number of drones required to meet either an average or 90th percentile response time goal while accounting for an $M/M/d$ queue using a chance constraint that guarantees a sufficient number of drones are located at each base. We leverage a conditional value-at-risk (CVaR) framework to control the tail of the response time distribution and we extend previous research by explicitly accounting for baseline response times (Section~\ref{Red}).

\item We develop a reformulation technique that exploits the baseline response times to induce coverage-like sparsity for median-type problems, allowing us to solve large-scale, real-world instances to optimality using an off-the-shelf solver (Section~\ref{Reform}). We use a two-stage solution pipeline that allows us to determine the minimum number of required drone resources to satisfy the constraints and assign those drones to OHCAs to minimize response times (Section~\ref{Imp}).

\item We demonstrate the application of our framework to determine the optimal deployment of AED-enabled drone resources to meet various response time targets using eight years of real data from a region covering 26,000 square kilometers around Toronto, Canada. These problems, which we solve to optimality, represent some of the largest instances in the literature for facility location models with queuing constraints. To generate the demand parameters and test our optimization model, we develop a two-stage machine learning approach to simulate cardiac arrest incidents. We first estimate the spatial OHCA distribution using kernel density estimation (KDE). Conditional on the location of an OHCA, we use a modified $K$-nearest neighbors algorithm to estimate its corresponding baseline response time (Section~\ref{NumExp}). Using our simulation-optimization framework, we show that:

\begin{itemize}
\item A modest number of drone resources is sufficient to significantly reduce response times in all regions. An objective function that minimizes average response time results in drone resources concentrated in cities, with little impact on the tail of the distribution. In contrast, optimizing for the tail of the response time distribution produces larger and more geographically dispersed drone networks that also manage to significantly reduce the average response time. 

\item Except for the most urban region, the total number of drone bases is not sensitive to changes in call volume. Additionally, the total number of drone resources in the large rural regions is insensitive to call volume. 

\item The optimal drone networks are insensitive to an objective function that models the trade-off between minimizing the total number of drone bases and minimizing the total number of drones. As a result, EMS providers may not need to consider the cost reductions associated with centralizing drone resources when designing their networks.
\end{itemize}

\end{enumerate}

Although we have described the problem in the context of AED-drones, both our model and reformulation technique are generalizable to any location-queuing problem where transportation options are added on top of an existing system to optimize performance. Other example applications include ambulance location, health clinic location, logistics problems with multiple transportation options (e.g., air vs. road vs. train), and delivery problems faced by organizations such as UPS and FedEx. 


\section{Literature Review}\label{Lit}

This paper contributes to the literature on facility location with congestion, the facility location literature applied to EMS, the medical literature on cardiac arrest response optimization, and the emerging operations research literature on drone applications.

Our work primarily contributes to the rich literature on facility location problems with congestion. Initial research in this domain focused on tractably integrating probabilistic queuing constraints within the classical integer programming framework for facility location, which can be partitioned into two groups depending on whether the facility location model employs a coverage or median-style framework. For coverage models, the pioneering work by \cite{ReVelle1988, ReVelle1989a, ReVelle1989b} incorporated congestion by assuming that the probability of two servers being busy at the same facility are independent. This assumption was later relaxed by \cite{Marianov1994, Marianov1996} who proposed a formulation to model the behaviour in each facility's catchment region as $M/M/d$ or $M/G/d$ queuing systems. The first median model to include queueing constraints -- the stochastic queue median -- was proposed by \cite{Berman1985} and employed a $M/G/1$ queue at each facility. \cite{Batta1989} later extended the stochastic queue median model to include multiple servers using an $M/G/d$ queue. However, the aforementioned models only \emph{implicitly} consider system congestion. That is, the number of servers at each location is fixed a priori.

In the literature on facility location with congestion, our work is most similar to that of \cite{Marianov2002}, who propose a set cover model where each base constitutes an $M/M/d$ queue and the number of servers is a variable to be optimized, i.e., \emph{explicit} congestion. Congestion is modeled using a chance constraint that considers the probability of a certain service level being met. To solve their model, they develop a heuristic algorithm and apply it to a single instance with 55 nodes. We employ similar queuing constraints as \cite{Marianov2002}, but in a median-style model, which is much more challenging to solve. Explicit congestion has also been considered in median-style models by \cite{Berman2007} and \cite{Aboolian2008}. However, our approach differs from theirs in two key ways: 1) we model mobile servers that travel to the demand, while they focus on fixed server locations (i.e., ATMs) that customers travel to, and 2) we model queuing using chance constraints on the queue length, rather than constraints on the average waiting time. Median models are known to be more computationally challenging compared to coverage models due to the lack of sparsity in the response time matrix \citep{Li2011}. To address this challenge, we exploit the historical (baseline) response times in our emergency response application to induce sparsity, allowing us to optimally solve large-scale, real-world problem instances. To the best of our knowledge, previous research has not attempted to include or exploit historical response times in this way.

Facility location models have been applied extensively to EMS location problems with the majority of previous research focusing on ambulances. For a comprehensive review of facility location in the context of EMS, please see \cite{Li2011}, \cite{Basar2012}, or \cite{Ahmadi2017}. Our work intersects with two streams of the EMS location literature: response time optimization and ambulance congestion. Models for response time optimization include objective functions that optimize the median/average response time (e.g., \cite{Carson1990, Serra1998}), the tail of the response time distribution (e.g., \cite{Krishnan2016, Chan2018}), and survival directly (e.g.,  \cite{Erkut2007, Knight2012}). 
Ambulance congestion models typically leverage queuing theory and chance constraints to handle demand uncertainty (e.g., \cite{Restrepo2008, ToroDiaz2013, Zayas2013}). Most location-queuing models are intractable and solved using heuristic approaches \citep{Brotcorne2003, Ahmadi2017}. 
Our paper differs from previous EMS location research because we combine a response time optimization model with an explicit queuing model to handle system congestion. Moreover, we exploit historical response times to induce sparsity in a novel reformulation technique, allowing us to solve large-scale problems to optimality using an off-the-shelf solver, instead of using heuristics.

There is a growing body of work that aims to improve treatment for cardiac arrest by developing optimization-based approaches for locating static AEDs in public locations. Most of this work utilizes coverage-based models, motivated by an assumed maximum distance that a bystander can travel to retrieve the AED before emergency responders arrive at the scene \citep{Chan2013, Sun2016, Chan2016}. More recently, \cite{Chan2018} developed a response-time based model using conditional value-at-risk to optimize the tail of the response time distribution for static AEDs under cardiac arrest location uncertainty. However, that model is tailored for locating static AEDs and does not address the other two limitations mentioned previously, namely system congestion and the historical response time distribution, which are important when optimizing response to private location OHCAs.

Lastly, our work contributes to the growing operations research literature on drone applications. The majority of this literature has focused on the problem of routing one or more drones in a coordinated fashion \citep{Xia2017} or on the flying sidekick travelling salesman problem, where a drone is routed in coordination with a traditional delivery vehicle \citep{Murray2015, Ferrandez2016, Dorling2017, Carlsson2018, Agatz2018}. Although many potentially impactful healthcare applications of drones have been proposed \citep{Scott2017}, the literature on such drone applications is surprisingly sparse. To the best of our knowledge, there are only five published studies that consider healthcare applications for drones to date. Besides the three drone studies mentioned above, there have been two studies examining delivery of medical supplies to rural areas \citep{Scott2017, Kim2017}.

\section{Models}\label{Models}

In this section, we present a model that integrates response time optimization with queuing (Section~\ref{Red}) and describe how baseline response times can be exploited to improve the computational efficiency of our model (Section~\ref{Reform}). Finally, we outline a two-stage solution pipeline that determines the minimum number of drones required to satisfy a response time constraint and the optimal assignments for those drone resources (Section~\ref{Imp}).

\subsection{Integrating response time and queuing optimization}\label{Red}

We present a general optimization framework that combines the $p$-median model with an explicit $M/M/d$ queue, where $d$ is modeled using a decision variable. We use our framework to minimize the total number of drones required to meet either an average or 90th percentile response time goal while guaranteeing that a certain service level, quantified by drone availability when an emergency call comes in, is met. In particular, we model the service level using a chance constraint on the queue length. Our model considers baseline response times so our model strategically locates drone resources to optimize system performance as a whole, considering the other non-drone resources available to respond to emergencies. The inclusion of this baseline response time data leads to significantly improved model tractability. We focus our initial exposition on an average response time objective and later extend the model to optimize the 90th percentile using CVaR.

Let $\mathcal{I}$ index the set of $m$ candidate drone base locations and $\mathcal{J}$ index a set of $n$ demand points (OHCAs) to be served. We add an artificial node $i^B$ to the set $\mathcal{I}$ to represent all current stations from which emergency response arises (i.e., ambulance bases, fire stations). The response time from $i^B$ to each OHCA is then derived from historical response times. The parameter used to represent the response time from drone base $i$ to demand point $j$ is given by $r_{ij}$ and the parameter used to represent the baseline response time to demand node $j$ is denoted $b_j$ (i.e., $r_{i^Bj} = b_j$). Let $\mathcal{D}_i$ be an index set for the number of drones, $d$, deployed at base $i$; for $i^B, \mathcal{D}_{i^B} = \{1\}.$ We let $x_{ij}$ represent a continuous assignment variable that represents the fraction of demand at node $j$ that is assigned to a drone at base $i$. The location decision variable is denoted by $y_{id}$ where $y_{id} = 1$ indicates that at least $d_i$ drones are located at base $i$.

To model system congestion, we consider drone base $i$ as an independent queuing system with $d_i$ servers. Each base defines a catchment area according to its share of assigned demand points. We assume that a Poisson process with an arrival rate of $\lambda_i$ governs the occurrences assigned to base $i\in\mathcal{I}$ and that the drone service time is an exponentially distributed random variable with parameter $\mu$. In contrast to ambulances, the drone travel time represents a smaller component of the total service time \citep{Boutilier2017}, so the memoryless assumption is more realistic. Let the state space $S=\{0,1,2,...\}$ denote the number of calls in the system. Let $\rho = \lambda/d\mu$ (we drop the $i$ from $d_i$ and $\lambda_i$ for the remainder of this paragraph for simplicity) and, assuming $\rho < 1$, let $\pi_s$ denote the steady-state probability associated with state $s$, which can be determined from solving the well-known, steady-state equations of an $M/M/d$ queue (Kleinrock 1975):
$$\pi_0 = \left[\sum_{s=0}^{d-1}\frac{(d\rho)^s}{s!}+\frac{(d\rho)^d}{d!}*\frac{1}{1-\rho}\right]^{-1}$$
\[
    \pi_s =
\begin{cases}
    \frac{d^s\rho^s\pi_0}{s!}, & s=1,2,...,d-1,\\
    \frac{d^d\rho^s\pi_0}{d!}, & s=d,d+1,...
    \end{cases}
\]
Queuing is incorporated into our optimization model through a service level constraint enforcing the steady-state probability at least one drone is available when an emergency occurs to be at least $\psi$. This constraint can be written as
$$\pi_0+\sum_{s=1}^{d-1}\pi_s \geq \psi,$$
or equivalently \citep{Marianov1998}
\begin{equation}
\sum_{s=0}^{d-1}\frac{(d-s)d!}{s!}\frac{1}{\rho^{d-s}} \geq \frac{1}{1-\psi}.\label{eq:servlevel}
\end{equation}
Let $\rho^{\psi}_{id}$ equal the value of $\rho$ that achieves equality in equation~\eqref{eq:servlevel} for base $i$. We represent equation~\eqref{eq:servlevel} in our optimization problem using the inequality
\begin{equation}\label{eq:congestion}\
\sum_{j\in\mathcal{J}}f_j x_{ij} \leq \mu \left(y_{i1}\rho^{\psi}_{i1}+\sum_{d\in\mathcal{D}_i}y_{id}(\rho^{\psi}_{id} - \rho^{\psi}_{i(d-1)})\right), \quad  i\in\mathcal{I},
\end{equation}
where $f_j$ is a scaling factor that converts a raw number of demand points assigned to base $i$ into a daily arrival rate. Intuitively, the left-hand side of~\eqref{eq:congestion} represents the daily number of OHCAs assigned to base $i$, while the right-hand side models the daily number of OHCAs that can be served at base $i$. In particular, $\rho^{\psi}_{id}$ is defined as the server utilization required to satisfy~\eqref{eq:servlevel}  when $d$ drones are located at base $i$, meaning that the sum on the right-hand side represents the incremental server utilization provided from adding an additional drone resource to base $i$. Since $\rho^{\psi}_{id}=\lambda/d\mu$, the right-hand side of~\eqref{eq:congestion} reduces to a bound on the arrival rate for base $i$. 

We define $\gamma$ as the desired improvement in seconds over the average baseline response time $\frac{1}{|\mathcal{J}|}\sum_{j\in\mathcal{J}} b_j$. The complete model that optimizes AED-enabled drone response, considering queuing and baseline response times is
\allowdisplaybreaks
\begin{subequations} \label{Med_RT}
\begin{align}
\quad\quad\underset{\by,\bx}{\mathrm{minimize}} \hspace*{1em} & \sum_{i\in\mathcal{I}}\sum_{d\in\mathcal{D}_i} y_{id}\\
\mbox{subject to}\hspace*{1em} & \frac{1}{|\mathcal{J}|}\sum_{j\in\mathcal{J}}\sum_{i\in\mathcal{I}} r_{ij}x_{ij} \leq \frac{1}{|\mathcal{J}|}\sum_{j\in\mathcal{J}} b_j- \gamma, \label{Cstr_RT}\\
&  \sum_{i\in\mathcal{I}} x_{ij}=1, \quad j\in\mathcal{J} , \label{Cstr_c}\\
& x_{ij} \leq y_{i1}, \quad i\in\mathcal{I}, j\in\mathcal{J},\label{Cstr_d} \\
& y_{id} \leq y_{i(d-1)}, \quad  d \in \mathcal{D}_i, i\in\mathcal{I}, \label{Cstr_e}\\
& \sum_{j\in\mathcal{J}}f_j x_{ij} \leq \mu \left(y_{i1}\rho^{\psi}_{i1}+\sum_{d\in\mathcal{D}_i}y_{id}(\rho^{\psi}_{id} - \rho^{\psi}_{i(d-1)})\right), \quad  i\in\mathcal{I},\label{Cstr_f}\\
& y_{i^B1} =1, \label{Cstr_g}\\
& 0\leq x_{ij} \leq 1, \quad i\in\mathcal{I}, j\in\mathcal{J}, \\
& y_{id}\in\{0,1\}, \quad  d \in \mathcal{D}_i,i\in\mathcal{I}.
\end{align}
\end{subequations}

The objective minimizes the total number of drones. Constraint \eqref{Cstr_RT} guarantees that the expected response time of the optimized drone network improves upon the average baseline response time by at least $\gamma$ seconds. Constraint \eqref{Cstr_c} ensures that each demand point is fully assigned. Constraints \eqref{Cstr_d} and \eqref{Cstr_e} are logical constraints that ensure a demand point can only be assigned to an open base and that the number of drones located at a base is properly modeled by the variables $y_{id}$ (i.e., $d-1$ drones must be allocated to a base, before $d$ drones can be allocated), respectively. Constraint \eqref{Cstr_f} determines the number of drones required at each base to meet the service level $\psi$, which depends on both the arrival and service rates. Constraint \eqref{Cstr_g} ensures that existing emergency response stations are included in the model. The implication is that response times will be no worse than the baseline. For conciseness in the remainder of the paper, we define the feasible region for formulation~\eqref{Med_RT} excluding constraint~\eqref{Cstr_RT} using
$$\bY = \Biggl\{\by\in\{0,1\}^{m\times |\mathcal{D}_i|} \; \bigr\rvert \; y_{i^B1}=1; y_{id} \leq y_{i(d-1)}, d\in\mathcal{D}_i,i\in \mathcal{I} \Biggr\}$$
and
\begin{multline*}
\bX(\by) = \Biggl\{\bx\in[0,1]^{m\times n} \; \bigr\rvert \; \sum_{i\in\mathcal{I}} x_{ij}=1, j\in\mathcal{J}; \; x_{ij} \leq y_{i1}, i\in\mathcal{I}, j\in\mathcal{J}; \\ \sum_{j\in\mathcal{J}}f_j x_{ij} \leq \mu \left(y_{i1}\rho^{\psi}_{i1}+\sum_{d\in\mathcal{D}_i}y_{id}(\rho^{\psi}_{id} - \rho^{\psi}_{i(d-1)})\right), i\in\mathcal{I}\Biggr\}.
\end{multline*}

\cite{Marianov2002} were the first to propose the modeling approach shown in constraint~\eqref{Cstr_f}, in the context of a coverage-type location model. Even though coverage-type models enjoy significant sparsity, they asserted that their model was ``especially difficult'' to solve and developed a heuristic solution approach. Our model adds to the difficulty because it includes the response-time constraint \eqref{Cstr_RT}, which leads to a fully dense input matrix consisting of the drone-based response time between each $i-j$ pair. To address this computational challenge, we develop a novel reformulation technique, which we describe next.

\subsection{Exploiting baseline response times}\label{Reform}

In this subsection, we describe how baseline response times can be exploited to improve computational efficiency of our model. The key idea is that an OHCA will never be assigned to a drone base with a response time that is worse than its baseline response time. Thus, even though the response time matrix is fully dense in theory, practically, we only need to consider a small fraction of all $i-j$ pairs. Mathematically, we can express this idea as follows.
\begin{lemma}\label{Prop1}
For any $j \in \mathcal{J}$, if $b_j \leq r_{ij}$ for a given $i \in \mathcal{I}\setminus\{i^B\}$, then there exists an optimal solution to~\eqref{Med_RT} such that $x^*_{ij}=0$.
\end{lemma}
\proof{Proof.}
Let $(\by^*,\bx^*)$ be an optimal solution to~\eqref{Med_RT}. Without loss of generality, assume that for a particular $\hat{j}\in\mathcal{J}$ and $\hat{i}\in\mathcal{I}\setminus\{i^B\}$, that  $x^*_{\hat{i}, \hat{j}}=1$ and $b_{\hat{j}} \leq r_{\hat{i}\hat{j}}$. Since, $b_{\hat{j}}=r_{i^B\hat{j}} \leq r_{\hat{i}\hat{j}}$, we can set $x^*_{i^B\hat{j}}=1$ and $x^*_{\hat{i}, \hat{j}}=0$, while maintaining the feasibility of constraint~\eqref{Cstr_RT}. Constraint~\eqref{Cstr_c} remains feasible because $\hat{j}$ is assigned to $i^B$. The remaining constraints are not affected and the $\by^*$ variable is unchanged, preserving the optimal cost.\halmos 
\endproof
Furthermore, if no drone base can improve on the historical response time, then the OHCA will be assigned to the baseline node $i^B$. Practically, this situation is more likely to occur in dense urban areas where ambulance response is faster.
\begin{corollary}\label{Prop2}
If $b_j \leq r_{ij}$ for all $i \in \mathcal{I}$, then there exists an optimal solution such that $x^*_{i^Bj} = 1$ and $x^*_{ij}=0, \forall i\neq i^B$.
\end{corollary}
The implicit assumption in this result is that our baseline node has infinite capacity. If the baseline is accurately derived from real data, the capacity of the system and associated congestion will already be accounted for in the historical response times, so this assumption will be reasonable. Another way to justify this assumption is that OHCA is the highest priority of all emergency calls. As a result, there will always be capacity in the system to respond to such calls, and even ambulances en route to another emergency will be re-routed to an OHCA if needed.

The sparsity that Lemma~\ref{Prop1} and Corollary~\ref{Prop2} allow us to induce can be implemented in a couple ways. One approach is to add additional constraints on $\bx$ in the model. However, this approach may be limited by that fact that the response time parameter matrix remains completely dense. Although modern solvers typically exploit this information during the pre-processing phase of setting up the problem, we found that manually pre-processing the parameters to remove all $i-j$ pairs that will eventually be set to zero results in a speed-up (see Results section). We believe this speed-up is due primarily to the fact that we avoid passing a fully dense matrix to the solver.

A second approach involves manipulating the parameter matrix and reformulating~\eqref{Med_RT} to maximize response time \emph{improvement}. To do this, we define $t_{ij} = \max\{b_j-r_{ij},0\},\forall i\in\mathcal{I}, j\in\mathcal{J}$ as the response time improvement over the baseline. This reformulation induces a significant amount of sparsity regardless of pre-processing because any $i-j$ pair that does not improve upon the baseline response time is assigned a response time improvement of zero. In addition, all assignments to the baseline service ($i^B$) receive a response time improvement of zero, which is in contrast to the previous approaches, where all $i^B-j$ pairs are represented by a nonzero response time. The transformed model, which we call the \emph{response time improvement model}, is written
\allowdisplaybreaks
\begin{subequations} \label{IRT}
\begin{align}
\quad\quad\underset{\by,\bx}{\mathrm{minimize}}  \hspace*{1em} & \sum_{i\in\mathcal{I}}\sum_{d\in\mathcal{D}_i} y_{id}\\
\mbox{subject to}\hspace*{1em} &  \frac{1}{|\mathcal{J}|}\sum_{j\in\mathcal{J}}\sum_{i\in\mathcal{I}} t_{ij}x_{ij} \geq \gamma, \label{Cstr_IRT}\\
&  \bx\in\bX(\by), \by \in \bY.
\end{align}
\end{subequations}

Constraint~\eqref{Cstr_RT} is replaced by constraint~\eqref{Cstr_IRT}, which guarantees that the average response time improvement is at least $\gamma$. The rest of the constraints are unchanged. Since we're dealing with average response time, it should be intuitive that formulation~\eqref{IRT} is equivalent to formulation~\eqref{Med_RT}. We state this equivalence formally below, after defining some notation, which also allows us to derive related results beyond average response time.

Let $\xi\in\Omega$ be a random vector representing the location of the next cardiac arrest event. With a slight abuse of notation, let $\by_B(\xi)$ denote the closest baseline emergency response base location to $\xi$ and $\by_{R\cup B}(\xi)$ denote the closest base to $\xi$ among both the drone and baseline locations. Using this notation, we define $R(\xi, \by_{R\cup B}(\xi))$ as the response time between the next demand arrival $\xi$ and the nearest combined (drone plus baseline) network location, $B(\xi,\by_B(\xi))$ as the response time between the next demand arrival $\xi$ and its nearest baseline location, and $I(\xi, \by_{R\cup B}(\xi))$ as the improvement in response time between $\xi$ and its nearest combined network location. Note that $I(\xi, \by_{R \cup B}(\xi)) = B(\xi, \by_B(\xi)) - R(\xi, \by_{R\cup B}(\xi))$ and by definition $B(\xi, \by_B(\xi)) \geq R(\xi, \by_{R\cup B}(\xi))$. Since, $R, B$, and $I$ are random variables, there are associated distributions of response time induced by the combined network locations, response time induced by the baseline locations, and response time improvement induced by the combined network locations, respectively. The expected values of $R$, $B$, and $I$ are given by $\mathbb{E}_{\xi}(R)$, $\mathbb{E}_{\xi}(B)$, and $\mathbb{E}_{\xi}(I)$, respectively.

The equivalence between constraints~\eqref{Cstr_RT} and~\eqref{Cstr_IRT} is straightforward to establish due to linearity of the expectation operator (proof omitted).

\begin{lemma}\label{Lem2}
$\mathbb{E}_{\xi}(I) \ge \gamma$ if and only if $\mathbb{E}_{\xi}(R) \leq \mathbb{E}_{\xi}(B) - \gamma$.
\end{lemma}

With these constraints being the only difference between formulation~\eqref{Med_RT} and~\eqref{IRT}, it follows that the two formulations are equivalent (proof omitted).

\begin{theorem}\label{Thm1}
A solution is optimal for~\eqref{Med_RT} if and only if it is optimal for~\eqref{IRT}.
\end{theorem}

Beyond average response time, we can also establish a related result for conditional value-at-risk (CVaR), which is a commonly used approximation for tail value-at-risk (VaR) metrics. The upper CVaR for a specified probability level $\beta$ in $(0,1)$ is defined as
$$\phi^U_{\beta}(R) = \min_{\alpha_{\beta}\in\mathbb{R}}\;\left\{ \alpha_{\beta} + \frac{1}{(1-\beta)}\mathbb{E}_{\xi}[R - \alpha_{\beta}]^+\right\},$$
where $[\cdot]^+$ represents $\max\{\cdot,0\}$ \citep{Rockafellar2000, Rockafellar2002}. The response time model equivalent to~\eqref{Med_RT} but with a CVaR constraint on response time (instead of average response time) can be written as
\allowdisplaybreaks
\begin{subequations} \label{RTCVAR}
\begin{align}
\quad\quad\underset{\by,\bx}{\mathrm{minimize}}  \hspace*{1em} &  \sum_{i\in\mathcal{I}}\sum_{d\in\mathcal{D}_i} y_{id}\\
\mbox{subject to}\hspace*{1em} & \phi^U_{1-\beta}(R) \leq \phi^U_{1-\beta}(B)- \gamma,\label{Cstr_CVARRT} \\
&  \bx\in\bX(\by), \by \in \bY,
\end{align}
\end{subequations}
and, the \emph{response time improvement model} with a CVaR constraint can be written as
\allowdisplaybreaks
\begin{subequations} \label{ITCVAR}
\begin{align}
\quad\quad\underset{\by,\bx}{\mathrm{minimize}}  \hspace*{1em} &  \sum_{i\in\mathcal{I}}\sum_{d\in\mathcal{D}_i} y_{id}\\
\mbox{subject to}\hspace*{1em} & \phi^U_{\beta}(I)\geq \gamma, \label{Cstr_CVARIT}\\
&  \bx\in\bX(\by), \by \in \bY.
\end{align}
\end{subequations}
Unfortunately, in the CVaR case, there is no analogous equivalence between constraints~\eqref{Cstr_CVARIT} and~\eqref{Cstr_CVARRT}. However, it is possible to show that satisfying constraint~\eqref{Cstr_CVARIT} implies constraint~\eqref{Cstr_CVARRT} is satisfied.
\begin{lemma}\label{Lem3}
For $\beta\leq 0.5$, if $\phi^U_{\beta}(I)\geq \gamma$ then $\phi^U_{1-\beta}(R) \leq \phi^U_{1-\beta}(B)- \gamma.$
\end{lemma}
\proof{Proof.}
\begin{align*}
\gamma & \leq \phi^U_{\beta}(I) \\
& = \phi^U_{\beta}(B-R) \\
& \leq \phi^U_{\beta}(B) + \phi^U_{\beta}(-R) \\
& = \phi^U_{\beta}(B) - \phi^U_{1-\beta}(R) \\
& \leq \phi^U_{1-\beta}(B) - \phi^U_{1-\beta}(R),
\end{align*}
where the first inequality comes from subadditivity of CVaR and the second inequality comes from the fact that $\phi_\beta^U(\cdot)$ is increasing in $\beta$. \halmos
\endproof
If $\beta = 0$, then Lemma~\ref{Lem3} is equivalent to Lemma~\ref{Lem2}. Unfortunately, for nonzero $\beta$, there is a gap due to the inequalities, which leads to a weakened one-way implication. Although this result implies that the CVaR improvement model is not equivalent to the CVaR response time model, formulation \eqref{ITCVAR} can be used to provide an upper bound on the optimal value of formulation~\eqref{RTCVAR} because of the one-way implication (proof omitted).
\begin{theorem}\label{Thm2}
The optimal value for~\eqref{ITCVAR} provides an upper bound on the optimal value for~\eqref{RTCVAR}.
\end{theorem}

In our computational results, we solve the discretized version of~\eqref{RTCVAR} shown below
\allowdisplaybreaks
\begin{equation} \label{RTCVARD}
\begin{aligned}
\quad\quad\underset{\by,\bx,\bz, \alpha}{\mathrm{minimize}}  \hspace*{1em} &  \sum_{i\in\mathcal{I}}\sum_{d\in\mathcal{D}_i} y_{id}\\
\mbox{subject to}\hspace*{1em} &  \alpha + \frac{1}{(1-\beta)|\mathcal{J}|}\sum_{j\in\mathcal{J}}z_j \leq \gamma_{90}, \\
& z_j \geq \sum_{i\in\mathcal{I}}\sum_{j\in\mathcal{J}} r_{ij}x_{ij} - \alpha, \quad j\in\mathcal{J}, \\
& z_j \geq 0, \quad j\in\mathcal{J}, \\
&  \bx\in\bX(\by), \by \in \bY,
\end{aligned}
\end{equation}
where $\gamma_{90}=\phi^U_{1-\beta}(B)- \gamma$ represents the desired response time goal for the 90th percentile. To solve this model, we employ the manual pre-processing step described previously and add the constraints on $\bx$ to exploit the baseline response times, as outlined in Lemma~\ref{Prop1} and Corollary~\ref{Prop2}.

\subsection{Complete solution pipeline}\label{Imp}

In this subsection, we describe a two-stage solution pipeline that we employ to 1) determine the minimum number of required drone resources to satisfy the constraints, and 2) assign those drones to OHCAs to minimize response times. We use a two-stage approach because although formulation~\eqref{Med_RT} (or~\eqref{IRT}) will minimize the number drones needed to reduce the average response time, an optimal solution may not make the optimal $i-j$ assignments with respect to response times, since the latter is a constraint, not an objective. To overcome this issue, we first solve formulation~\eqref{Med_RT} (or~\eqref{IRT}) to determine the minimum number of drone resources, $P^*$, needed to satisfy the associated response time constraint. We then solve the second stage model below to determine optimized response times given $P^*$ drone resources:
\allowdisplaybreaks
\begin{equation} \label{IRTmax}
\begin{aligned}
\quad\quad\underset{\by,\bx}{\mathrm{maximize}} \hspace*{1em} &   \frac{1}{|\mathcal{J}|}\sum_{j\in\mathcal{J}}\sum_{i\in\mathcal{I}} t_{ij}x_{ij} \\
\mbox{subject to}\hspace*{1em} &  \sum_{i\in\mathcal{I}}\sum_{d\in\mathcal{D}_i} y_{id} = P^*, \\
&  \bx\in\bX(\by), \by \in \bY.
\end{aligned}
\end{equation}
The two-stage pipeline for the CVaR case is similar. We first solve formulation~\eqref{RTCVAR} (or~(\ref{RTCVARD})). We then solve the corresponding second stage model below:
\allowdisplaybreaks
\begin{equation} \label{RTCVARmax}
\begin{aligned}
\quad\quad\underset{\by,\bx,\bz,\alpha}{\mathrm{minimize}} \hspace*{1em} &    \alpha + \frac{1}{(1-\beta)|\mathcal{J}|}\sum_{j\in\mathcal{J}}z_j\\
\mbox{subject to}\hspace*{1em} &  z_j \geq \sum_{i\in\mathcal{I}}\sum_{j\in\mathcal{J}} r_{ij}x_{ij} - \alpha, \quad j\in\mathcal{J}, \\
& z_j \geq 0, \quad j\in\mathcal{J}, \\
& \sum_{i\in\mathcal{I}}\sum_{d\in\mathcal{D}_i} y_{id} = P^*, \\
&  \bx\in\bX(\by), \by \in \bY.
\end{aligned}
\end{equation}

Both first-stage and second-stage models can be solved to optimality in our large-scale instances due to the tractability induced by Lemma~\ref{Prop1} and Corollary~\ref{Prop2}. Without the inclusion of baseline response times, second-stage models~\eqref{IRTmax} and~\eqref{RTCVARmax} are equivalent to the $p$-median model with probabilistic queuing constraints, which is difficult to solve for large instances.

\subsection{Drones versus drone bases}

We further extend our models to investigate the trade-off between minimizing the number of drones versus minimizing the number of drone bases. Although currently unknown, the costs associated with building a new base versus adding a drone to a current base will likely be different. For example, EMS providers may prefer to centralize drones in fewer bases to simplify maintenance requirements and reduce fixed costs. In this case, minimizing the number of drone bases would be more appropriate. Alternatively, technological advances such as DroneBox \citep{DroneBox}, which provides a fully autonomous housing box that can charge and service drones, may make variable costs the dominant factor. In this case, minimizing the number of drones would be the appropriate objective. To study this mathematical trade-off, we use the objective function
\begin{equation}\label{DroneBase}\zeta\sum_{i\in\mathcal{I}}\sum_{d\in\mathcal{D}_i} y_{id} + (1-\zeta) \sum_{i\in\mathcal{I}} y_{i1},\end{equation}
where $\zeta\in[0,1]$. Note that $\zeta = 1$ corresponds to the original objective function.

\section{Case study: setup}\label{NumExp}

In this section, we provide details about our dataset (Section~\ref{data}) and model parameters such as candidate base locations (Section~\ref{bases}) and drone specifications (Section~\ref{dronespecs}). We then outline our machine learning approaches to simulate OHCA locations (Section~\ref{MLlocs}) and baseline response times (Section~\ref{MLrt}). Lastly, we outline the setup for our numerical optimization experiments (Section~\ref{ExpSetup}).

\subsection{Historical cardiac arrests}\label{data}

The Toronto Regional RescuNET comprises eight regions in Southern Ontario, Canada (Toronto, Durham, Simcoe, Muskoka, Peel, Hamilton, Halton, and York) with a total population of 7.54 million in a total area of 26,364 km$^2$. Each region is served by an independent paramedic service and there is a tiered response to emergency calls, where OHCA is the highest priority. Both paramedics and fire fighters are dispatched to all suspected OHCAs. 

We obtained Universal Transverse Mercator (UTM) coordinates and the historical 911 response time for all confirmed non-traumatic OHCA episodes throughout RescuNET from January 1, 2006 to December 31, 2014. Data was obtained from the Rescu Epistry cardiac arrest database \citep{Morrison2008, Lin2011} with research ethics approval. The Rescu Epistry data only includes confirmed historical OHCAs, which underestimates the total number of suspected OHCAs to which a drone would be dispatched. Although many suspected OHCAs will not turn out to be actual cardiac arrests (e.g., a ``life status questionable'' or ``unconscious'' call may turn out to be an intoxicated individual sleeping), EMS systems must still respond to those cases at the highest priority. We believe the same will be true for drones. Thus, we multiplied the annual number of confirmed OHCAs in each region by five as an estimate of the volume of suspected OHCAs -- this number was chosen through consultation with clinical collaborators in several regions and estimates that they provided. We then tested the sensitivity of our results to multipliers both smaller and larger than five. One rationale for a larger multiplier is that if they turn out to be effective, drones may get dispatched to even more calls at the highest priority level, even if the call is clearly not a cardiac arrest, in anticipation of a possible cardiac arrest. Such examples include calls for motor vehicle accidents or chest pain, where the individuals are alive but may degrade into cardiac arrest, for which an AED would be useful. 

Table~\ref{Table1} provides a summary of the eight RescuNET regions. Toronto, the largest city in Canada, is the most urban region, while Muskoka, nicknamed ``Cottage Country", is the most rural region. The other six regions include a mix of both dense urban cities and rural areas. For example, Peel region includes the sixth (Mississauga) and ninth (Brampton) largest cities in Canada, and Caledon, a sparse rural area with a population density of 86 people per km$^2$. We present our numerical results using Toronto, Muskoka, and Peel as representative urban, rural, mixed regions, respectively.

\renewcommand{\arraystretch}{1.5}
\begin{table}[h!]
\caption{Summary statistics and a comparison between the historical and simulated 911 response times for the eight regions comprising RescuNET.}\label{Table1}
\centering
\resizebox{\textwidth}{!}{%
\begin{tabular}{c c | c c c c c c c c }
\noalign{\smallskip} \hline
\multicolumn{2}{c|}{\multirow{2}{*}{\textbf{Characteristics}}} &\multicolumn{8}{c}{\textbf{Region}}\\ \cline{3-10}
&& Toronto & Peel & Simcoe & York & Halton & Durham & Hamilton & Muskoka \\ \hline
\multicolumn{2}{c|}{Population} & 2,731,571 & 1,381,739 & 479,650 & 1,109,909 & 548,435 & 645,862  & 536,917 & 103,423   \\
\multicolumn{2}{c|}{Population density (per km$^2$) }& 4334.4 & 1108.1	& 98.7& 629.9& 568.9	&	 255.9&480.6	&	7.8	\\
\multicolumn{2}{c|}{Average annual number of confirmed cardiac arrests }& 2977&848&440	&666&	355		&570 &618	&73  \\
\multicolumn{2}{c|}{Number of paramedic, fire, and police stations} &158&	68&	76&68	&41	&44 &	51	&32	 \\  \hline
\multirow{2}{*}{Historical 911 response time (s)} & Mean & 420 & 357  & 478 & 420 &  390& 364 & 400   &602  \\
& 90th percentile & 632 & 497& 780& 615 &  570& 540 & 660  & 1125\\
\multirow{2}{*}{Simulated training set 911 response time (s)} & Mean& 469 & 394 & 492 & 463 & 421 &421& 401  &655  \\
& 90th percentile & 625    &  509           & 807   & 595         &    562     &   514 & 481            &  1140   \\
\multirow{2}{*}{Simulated average of 100 testing set 911 response times (s)} & Mean & 469 & 395  &  493 &467   &  416   & 412 & 397   & 614  \\
& 90th percentile & 618     & 499     &       760   & 607         &  573    &  523  &   499            & 1149     \\ \hline
\multicolumn{2}{c|}{Optimal KDE bandwidth parameter ($b$) } & 112 &487  & 841 &565 & 766 &848 &826  &1412  \\
\multicolumn{2}{c|}{Optimal KNN hyper-parameters $(k,a,b)$ } & (380, 50, 0.65) &(160, 39, 0.90)  & (40, 27, 0.82) & (60, 48, 0.71) & (50, 31, 0.80)& (40, 54, 0.94)& (70, 0, 0.60)  & (20, 48, 0.86)   \\ \hline
\end{tabular}}
\end{table}

\subsection{Candidate base locations}\label{bases}

All paramedic, fire, and police stations within RescuNET were considered as candidate drone base locations. Addresses for each station were obtained from the regional provider and converted to UTM coordinates. The total number of candidate bases in each region is given in Table~\ref{Table1}.

\subsection{Drone parameters}\label{dronespecs}

Drone parameters used in our model were based on specifications reflecting recent
technological capabilities. Vertical acceleration/deceleration was set to 9.81 m/s$^2$ and horizontal acceleration/deceleration was set to 19.6 m/s$^2$ \citep{Schollig2011, Kumar2012}. Note that typically horizontal acceleration/deceleration is done simultaneously with vertical acceleration/deceleration. Maximum forward velocity was set at 27.8 m/s and the flying height was assumed to be 60m, which is below the maximum height allowed in Canada and the United States \citep{FAA,TransportCanada}. Accounting for the time it takes to reach maximum speed and height in an idealized situation, we assume 10 seconds is required for takeoff and landing.

\subsection{Estimating the spatial OHCA distribution}\label{MLlocs}

We estimate the spatial OHCA distribution separately for each region using kernel density estimation (KDE). KDE is a popular non-parametric approach for estimating the probability density function of a random variable from a finite sample \citep{Sheather1991}. It requires the specification of a kernel function and a bandwidth parameter, $b$. The bandwidth parameter controls the degree of smoothing, where a larger bandwidth value results in a smoother distribution. KDE has previously been used to estimate spatial OHCA distributions, which is supported by the fact that the spatial OHCA distribution has been shown to exhibit temporal stability \citep{Chan2016}.

We use $10$-fold cross validation to determine the optimal bandwidth parameter using Gaussian kernels for each region. We evaluate model performance using log-likelihood, which is a commonly used technique for both KDE \citep{Jones2009} and general statistical inference \citep{Bishop2006}. The optimal bandwidth parameters for each region are given in Table~\ref{Table1}. We use the optimal bandwidth parameter and all historical OHCA data in a region to fit a final KDE model for each region. We use this model to simulate a full year of confirmed OHCAs to serve as a training set and we simulate 100 different one-year test sets. Because of our assumption that drones will respond to all suspected cardiac arrests, we modify the arrival rate so that the full year of confirmed OHCAs occur during a 73 day period, representing a suspected OHCA multiplier of five. Figure~\ref{KDE} displays the study region, candidate base locations, and simulated training data.

\subsection{Estimating the 911 response time distribution}\label{MLrt}

In addition to simulating the location of OHCAs, we also estimate the 911 response time for each simulated OHCA. To do this, we use a modified $k$-nearest neighbors (KNN) algorithm. KNN is a non-parametric method for both classification and regression \citep{Altman1992}. In KNN-regression, the dependent variable is estimated using an average of the dependent variable value of the $k$ nearest neighbors. In our approach, we calculate a weighted average using the inverse of each neighbors distance to the target so that closer neighbors will have greater influence on the predicted outcome.

\begin{figure}[t]
\begin{center}
\subfigure[\  \label{KDE}]{
\frame{\includegraphics[width=.48\textwidth]{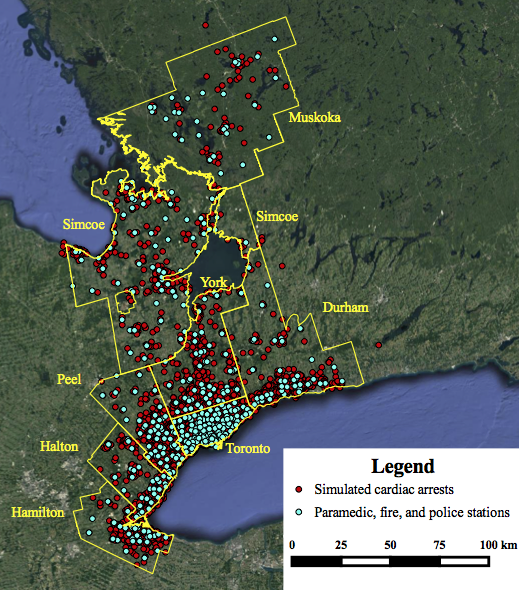}}}
\subfigure[\  \label{KNN}]{
\frame{\includegraphics[width=.4\textwidth]{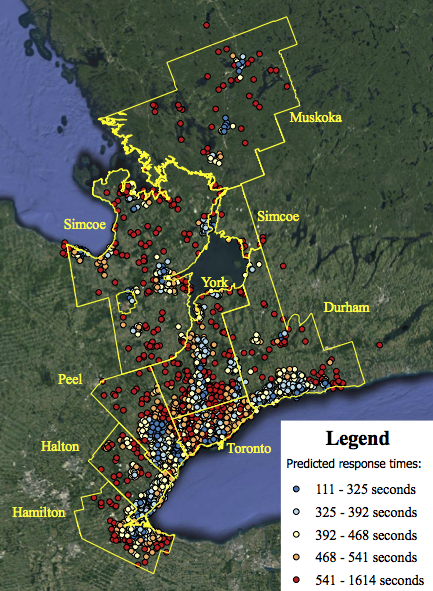}}}
\end{center}
\caption{(a) A map of the simulated cardiac arrests with all paramedic, fire, and police stations, and (b) A map of the simulated cardiac arrests and the estimated 911 response times. \label{MapPop}}
\end{figure}

We modify the traditional KNN algorithm in two ways to better fit our application. First, estimating 911 response times can be viewed as a bi-objective prediction problem because EMS agencies are typically measured on both mean and tail response time performance. Thus, we aim to minimize error in both the mean and upper tail of the estimated distribution. To do this, we create a simple scoring function equal to the sum of the mean absolute error and the absolute error at the 90th percentile. The 90th percentile is chosen because it is an internationally accepted metric for OHCA tail response times. Second, we apply a transformation that allows us more fine-grained control of the predicted response time distribution. In particular, let $\hat{\bh}$ represent a vector of historical 911 response times and let $\bar{\bh}$ represent the KNN-estimated 911 response times. The transformed 911 response time distribution $\bh$ is given by
$$ \bh = \mathbb{E}(\hat{\bh}) + a + (\bar{\bh} - \mathbb{E}(\hat{\bh}))\left(b \cdot\frac{\sigma(\bar{\bh})}{\sigma(\hat{\bh})}\right),$$
where $a$ (translational shift) and $b$ (degree of dispersion) are hyper-parameters that can be tuned and $\sigma$ represents the standard deviation. For each region, we use $10$-fold cross validation on the historical data to determine the optimal values of $k$, $a$, and $b$. Then, using the optimal hyper-parameters, we estimate the 911 response time for each simulated OHCA from Section~\ref{MLlocs}. Figure~\ref{KNN} displays the simulated training data color-coded according to the estimated 911 response time.

Although a modest contribution, we believe our simulation approach is valuable for several reasons. First, we do not have detailed historical data for all suspected OHCAs and as a result, we use this approach to obtain response time estimates for all simulated locations. Second, this approach can be used to reduce ``missingness'' in real data. For example, 8\% of our historical cardiac arrests were missing response times, and our approach provides a method to fill in those gaps. Third, we believe our general simulation framework will be useful for other researchers investigating spatial machine learning and optimization problems, because it can be used to generate large datasets for analysis and testing. Fourth, we use this approach to generate error bars by simulating many years of OHCA data in separate test sets.

\subsection{Experimental setup}\label{ExpSetup}

We focus our numerical experiments on two key metrics: the average response time and the 90th percentile response time. We treat each of the eight regions in RescuNET independently to mimic the current separation of EMS agencies. For the average response time, we solve the response time improvement model pipeline (formulations~\eqref{IRT} and~\eqref{IRTmax}), denoted \textbf{RTI-mean}, for values of $\gamma \in \{1,2,3\}$. For the 90th percentile response time, we solve the response time pipeline (formulations~\eqref{RTCVARD} and~\eqref{RTCVARmax}), denoted \textbf{RT-CVaR}, for 15\%, 30\%, and 50\% reductions in the 90th percentile response time. We use percentage reductions for the 90th percentile because of the large differences in the 90th percentile response time across the regions (see Table~\ref{Table1}). In all our models, we set the probability at least one drone is available when an emergency occurs, $\psi$, to 0.99. For both \textbf{IRT-mean} and \textbf{RT-CVaR}, we generate results to investigate the trade-off between minimizing the number of drone resources and the number of drone bases by using equation~\eqref{DroneBase} as the objective with $\zeta\in\{0,0.25,0.50,0.75,1\}$.

\section{Case study: results}\label{Results}

In this section, we present the numerical results from applying our framework to the problem of designing a drone network to deliver AEDs to OHCAs in Southern Ontario, Canada. In particular, we quantify the required drone network size to meet various response time improvement goals (Section~\ref{Network}), we conduct a sensitivity analysis on the drone network size as a function of the potential emergencies drones may be tasked with responding to (Section~\ref{Sens}), and we explore the trade-off between minimizing the total number of drones versus the total number of drone bases (Section~\ref{DBtrade}).

All optimization and machine learning experiments were programmed using \texttt{Python 3.5}. Optimization problems were solved using \texttt{Gurobi 7.0} with a maximum time limit of 10 hours and run on a desktop computer with an Intel Core i7-4790K 4.0 GHz processor and 32 GB of RAM. Optimal solutions for both \textbf{RTI-mean} and \textbf{RT-CVaR} were found within 1 hour for all instances, except Toronto, which required up to 10 hours for the \textbf{RT-CVaR} instances. For the average response time experiments across all eight regions, \textbf{RTI-mean} was able to reduce the solution time over a model that does not exploit baseline response times by 29.4\%, 13.0\%, and 17.7\% for $\gamma$ equal to 1, 2, and 3 minutes, respectively. For the 90th percentile CVaR experiments across all eight regions, \textbf{RT-CVaR} was able to reduce the solution time over a model that does not exploit baseline response timesby 25.4\%, 4.4\%, and 14.7\% for 90th percentile response time reductions of 15\%, 30\%, and 50\%, respectively.

\subsection{Drone network size}\label{Network}

Table~\ref{TableNetwork5} lists the number of drone bases and the total number of drones in each region prescribed by the \textbf{RTI-mean} and \textbf{RT-CVaR} models for each response time improvement goal. For example, to deliver an AED by drone 1 minute before 911 arrival on average, Toronto requires 1 drone base with 3 total drones, Muskoka requires 1 drone base with 1 drone, and Peel requires 1 drone base, with 3 total drones. As expected, the number of drone bases and total drones increases as the improvement goals become more ambitious. Furthermore, the number of drones required in the \textbf{RT-CVaR} solutions are significantly larger than in the \textbf{RTI-mean} solutions, since the goal is to reduce the tail of the distribution. For example, Durham requires 3 drone bases with 6 total drones to improve the average response time by 3 minutes, but even 44 drone bases (the maximum possible) is not able to reduce the 90th percentile by 50\%. Compared to a sequential location-queuing approach \citep{Boutilier2017}, our model is able to reduce the number of drone bases and total drones by 55\% and 63\%, respectively, for the same response time goal.

\begin{figure}[t]
\begin{center}
\subfigure[\ $\gamma=1$ \label{RTI1}]{
\frame{\includegraphics[width=.32\textwidth]{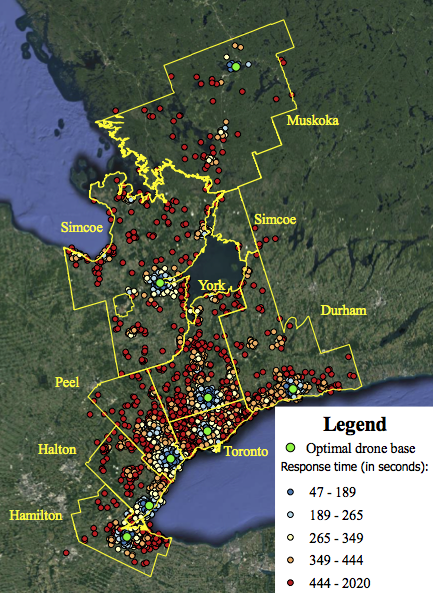}}}
\subfigure[\ $\gamma=2$ \label{RTI2}]{
\frame{\includegraphics[width=.32\textwidth]{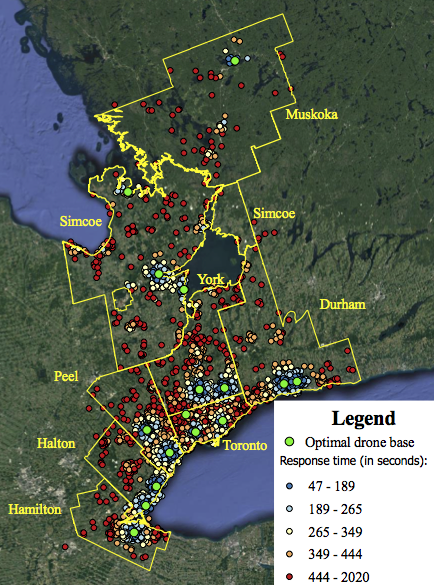}}}
\subfigure[\ $\gamma=3$ \label{RTI3}]{
\frame{\includegraphics[width=.32\textwidth]{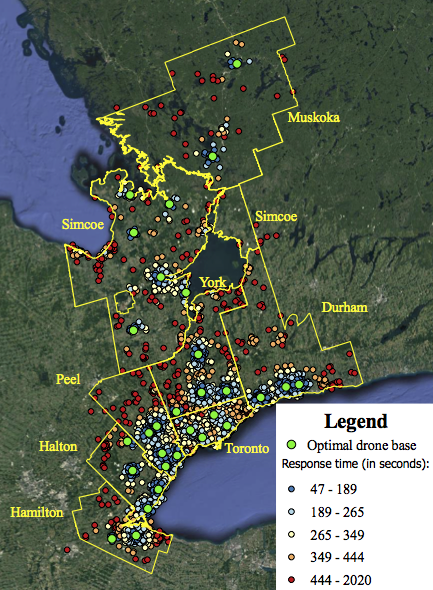}}}
\subfigure[\ 15\% \label{RT30}]{
\frame{\includegraphics[width=.32\textwidth]{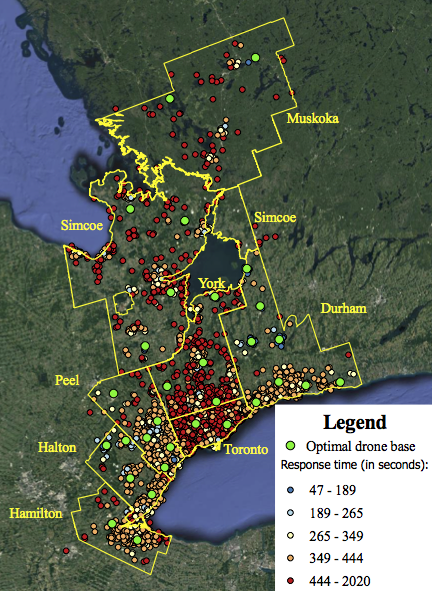}}}
\subfigure[\ 30\% \label{RT30}]{
\frame{\includegraphics[width=.32\textwidth]{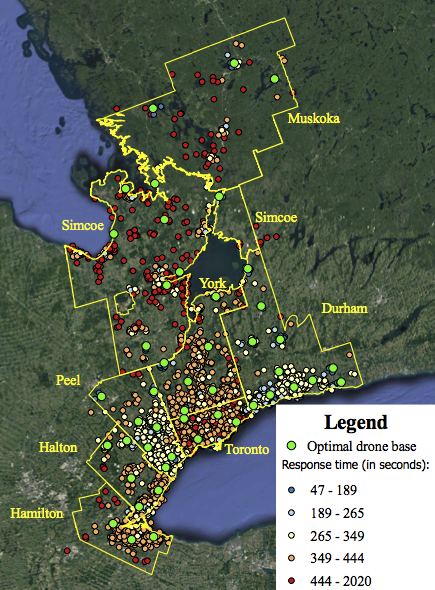}}}
\subfigure[\ 50\% \label{RT30}]{
\frame{\includegraphics[width=.32\textwidth]{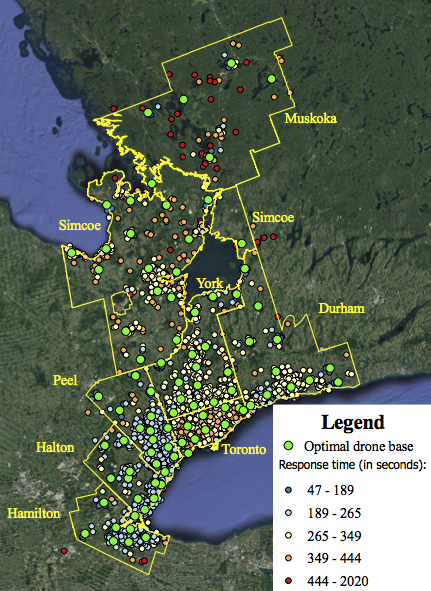}}}
\caption{Geographic layout of the optimized drone networks for \textbf{RTI-mean} (a-c) and \textbf{RT-CVaR} (d-f).\label{locs}}
\end{center}
\end{figure}

An important practical difference between \textbf{RTI-mean} and \textbf{RT-CVaR} is the distribution of drone resources. In particular, we find that \textbf{RTI-mean} tends to concentrate drone resources at fewer bases, while \textbf{RT-CVaR} distributes drone resources more widely across each region (see Table~\ref{TableNetwork5} and Figure~\ref{locs}). Furthermore, the regions with a larger geographic area tend to produce solutions with fewer drones per base because the models must distribute resources over a larger area to improve response times. For example, the smallest regions (e.g., Toronto and Hamilton) always require multiple drones per base, while the largest regions (e.g., Muskoka, Simcoe) only require one drone per base, unless the base is placed in a large city.

\renewcommand{\arraystretch}{1.25}
\begin{table}[h!]
\caption{Drone network size (Bases $\mid$ Total drones) for both \textbf{RTI-mean} and \textbf{RT-CVaR}.}\label{TableNetwork5}
\centering
\resizebox{\textwidth}{!}{%
\begin{tabular}{ c | c c c c c c c c}
\noalign{\smallskip} \hline
&\multicolumn{8}{c}{\textbf{Region}}\\ \cline{2-9}
 Improvement & Toronto & Peel & Simcoe & York & Halton & Durham & Hamilton & Muskoka   \\ \hline
&\multicolumn{8}{c}{\textbf{RTI-mean}} \\ \hline
\emph{$\gamma=1$} & 1 $\mid$ 3 & 1 $\mid$ 3 & 1 $\mid$ 2 & 1 $\mid$ 2 & 1 $\mid$ 2 & 1 $\mid$ 2 & 1 $\mid$ 2 & 1 $\mid$ 1\\
\emph{$\gamma=2$} & 3 $\mid$ 7 & 2 $\mid$ 5 & 3 $\mid$ 4 & 2 $\mid$ 4 & 2 $\mid$ 4 & 2 $\mid$ 4 & 1 $\mid$ 3 & 1 $\mid$ 1 \\
\emph{$\gamma=3$} & 4 $\mid$ 12 & 5 $\mid$ 10 & 6 $\mid$ 7 & 4 $\mid$ 7 & 4 $\mid$ 7 & 3 $\mid$ 6 & 2 $\mid$ 4 & 2 $\mid$ 2 \\
 \hline
 &\multicolumn{8}{c}{\textbf{RT-CVaR}} \\ \hline
15\% & 2 $\mid$ 4 & 5 $\mid$ 5 & 4 $\mid$ 4 & 4 $\mid$ 4 & 4 $\mid$ 4 & 8 $\mid$ 9 & 2 $\mid$ 3 & 2 $\mid$ 2\\
30\% & 4 $\mid$ 8 & 7 $\mid$ 9 & 7 $\mid$ 7 & 6 $\mid$ 7 & 5 $\mid$ 5 & 11 $\mid$ 13 & 3 $\mid$ 4 & 4 $\mid$ 4\\
50\% & 10 $\mid$ 16 & 14 $\mid$ 17 & 17 $\mid$ 17 & 15 $\mid$ 16 & 11 $\mid$ 11 & Infeasible & 7 $\mid$ 9 & 6 $\mid$ 6 \\
\hline
\end{tabular}}
\end{table}

Figure~\ref{locs} displays the geographical layout of the optimized drone networks for \textbf{RTI-mean} and \textbf{RT-CVaR} with all training set OHCAs colour-coded according to the optimized response time distribution. The locations determined by \textbf{RTI-mean} focus on areas with high OHCA incidence, corresponding to the largest cities in each region. By focusing on areas that include the majority of OHCAs, the model is able to efficiently reduce average response times. In contrast, the \textbf{RT-CVaR} base locations include both major cities and rural areas, which allows the network to impact the tail of the response time distribution. These differences highlight an equity-efficiency trade-off between the models, where \textbf{RT-CVaR} provides more equitable response time improvement across the regions and \textbf{RTI-mean} focuses on efficiently improving response times in major cities.

\begin{figure}[t]
\begin{center}
\subfigure[\ $\gamma=1$ \label{TestFigm1}]{
\frame{\includegraphics[width=.32\textwidth]{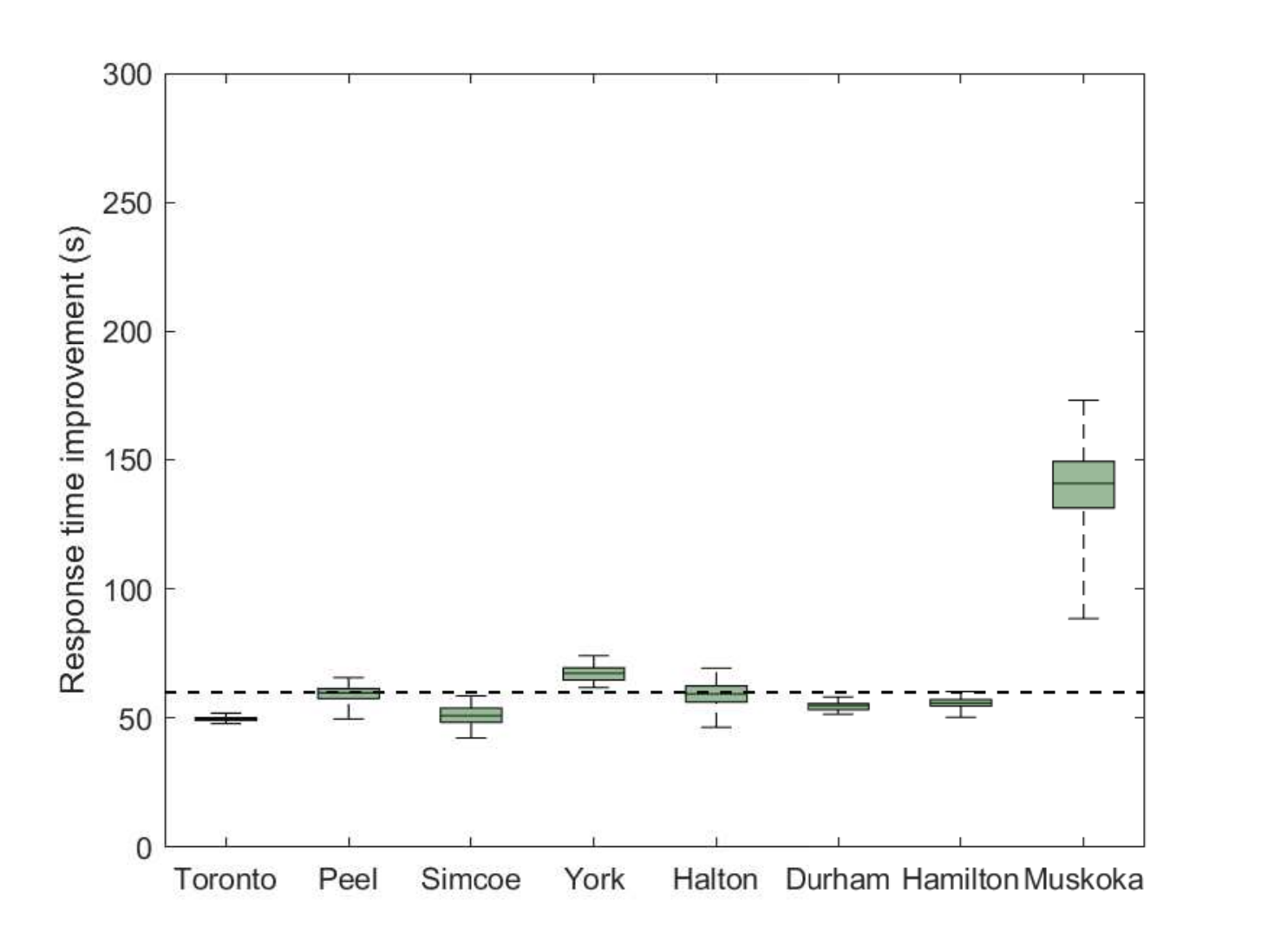}}}
\subfigure[\ $\gamma=2$ \label{TestFigm2}]{
\frame{\includegraphics[width=.32\textwidth]{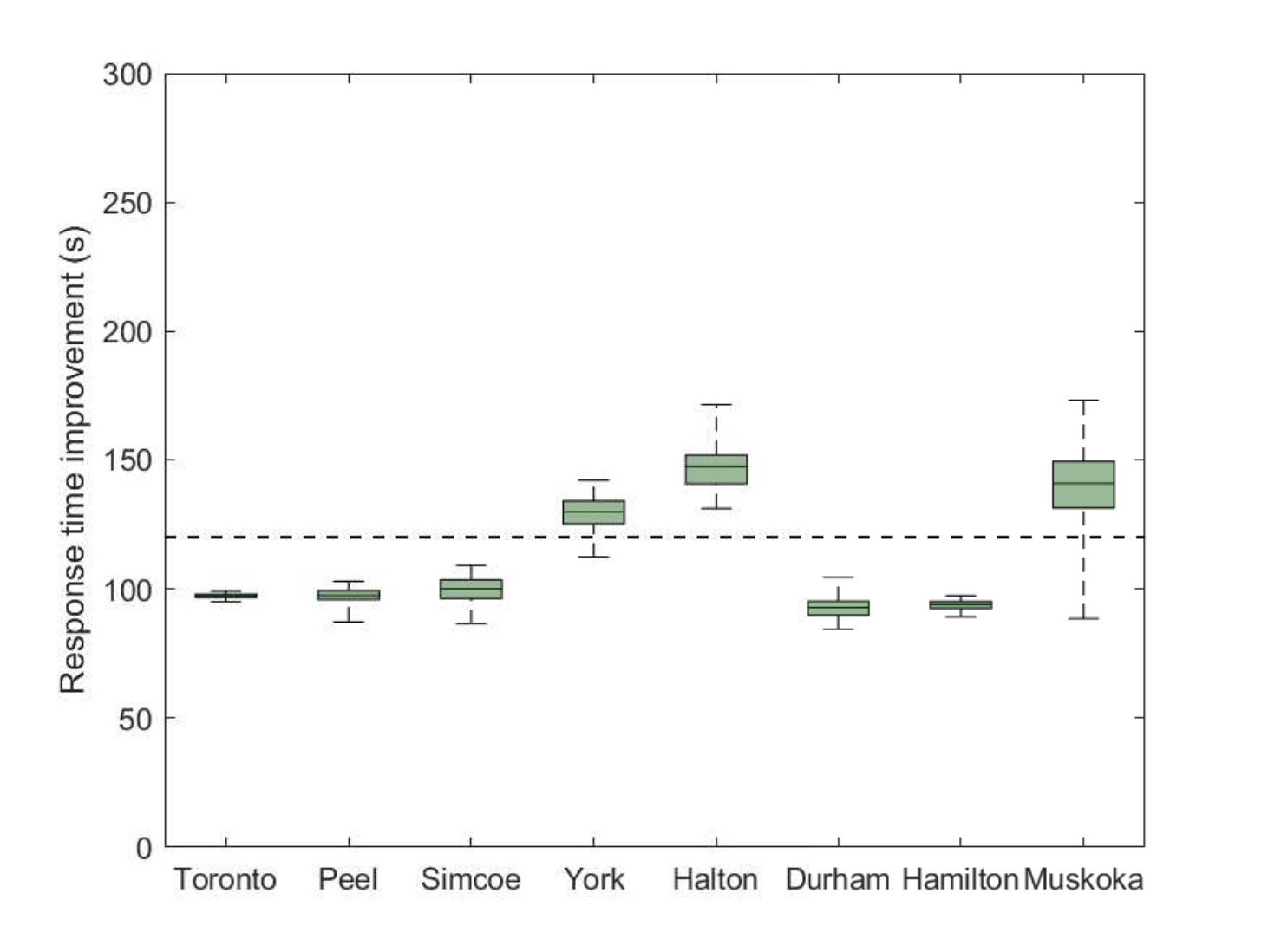}}}
\subfigure[\ $\gamma=3$ \label{TestFigm3}]{
\frame{\includegraphics[width=.32\textwidth]{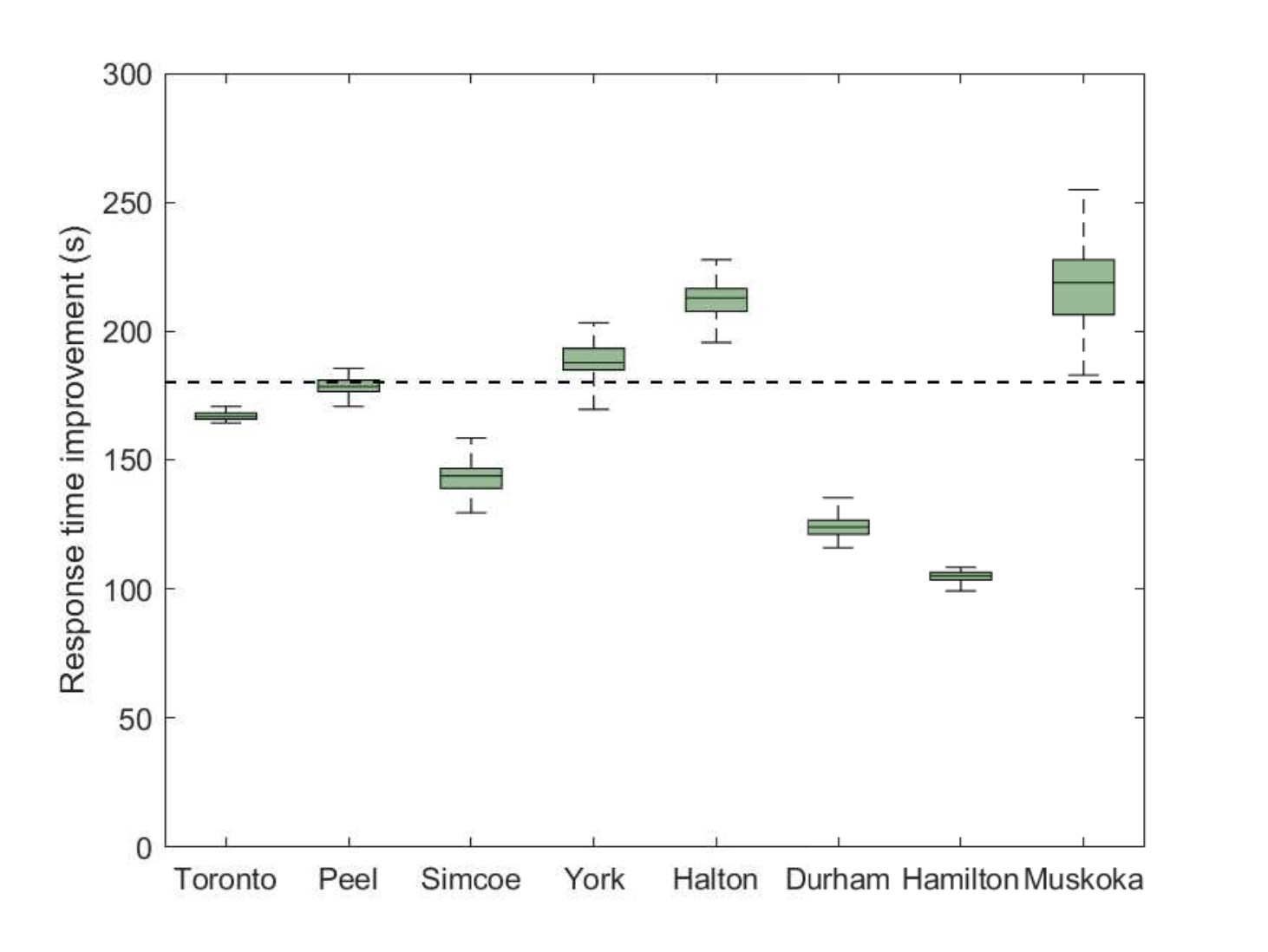}}}
\subfigure[\ 15\% \label{TestFigcv1}]{
\frame{\includegraphics[width=.32\textwidth]{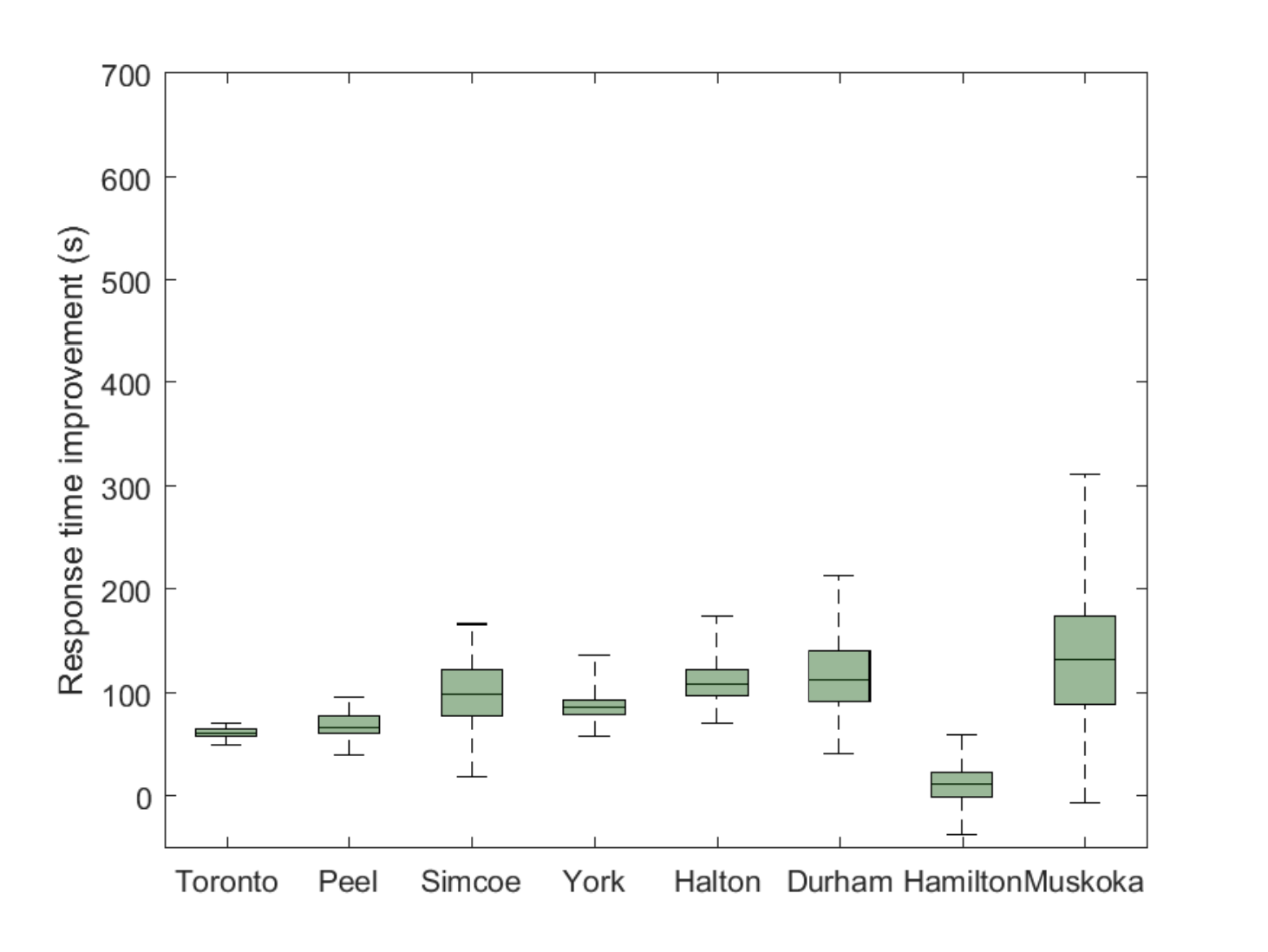}}}
\subfigure[\ 30\% \label{TestFigcv2}]{
\frame{\includegraphics[width=.32\textwidth]{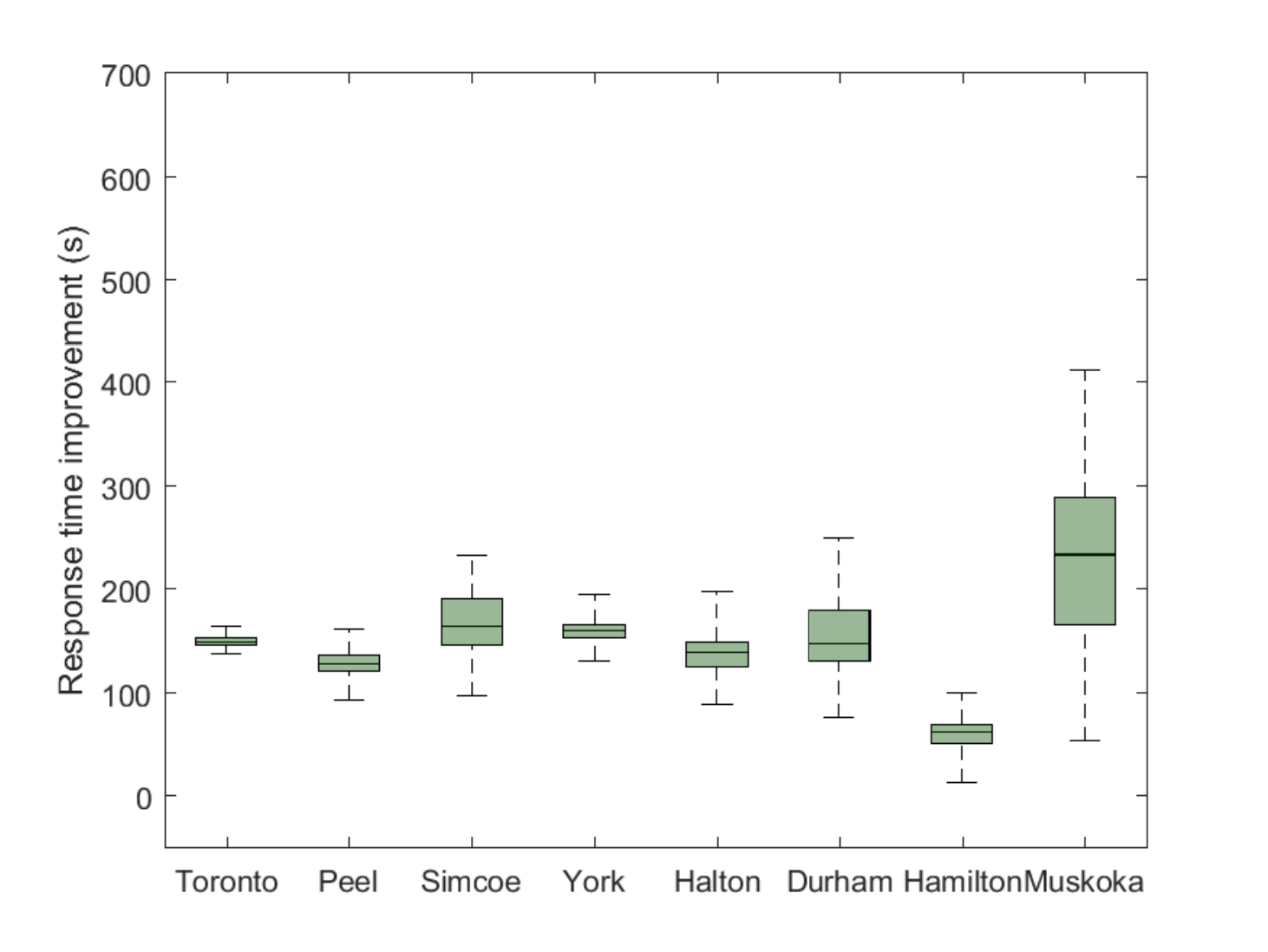}}}
\subfigure[\ 50\% \label{TestFigcv3}]{
\frame{\includegraphics[width=.32\textwidth]{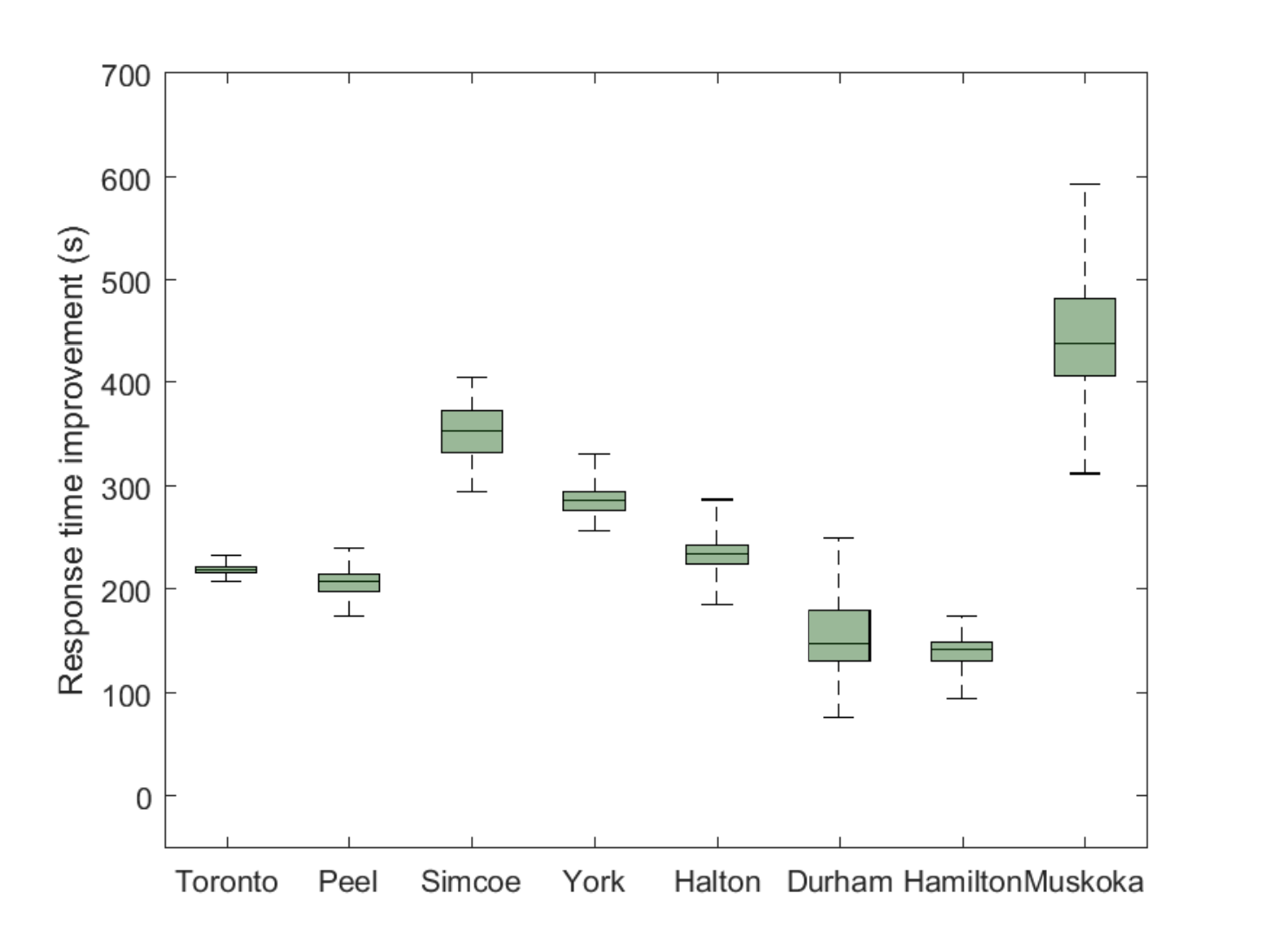}}}
\caption{Comparison of response time improvement for \textbf{RTI-mean} (a-c) and \textbf{RT-CVaR} (d-f) across 100 test sets. The dashed line in figures (a-c) represents the training set performance. \label{TestFig}}
\end{center}
\end{figure}

Figures~\ref{TestFig} displays the distribution of average and 90th percentile response time improvement across all 100 test sets. The test set improvement distributions are the same for Muskoka for $\gamma=1$ and $\gamma=2$ because the minimum drone network required to satisfy the queuing constraints for one minutes reduces the average response time by two minutes (Table~\ref{TableNetwork5} also shows that the networks are identical). The 50\% reduction in the 90th percentile in Durham is infeasible so we display the solution using the 30\% network. The range of test set performance for \textbf{RT-CVaR} (21s in Toronto to 359s in Muskoka) is much larger than \textbf{RTI-mean} (6s in Toronto to 92s in Muskoka).

For both models, the performance of the drone network in Muskoka is the most variable and we believe this is because of the small number of annual OHCAs. More specifically, geographical location differences in only a few OHCAs can significantly alter response time performance. Furthermore, the performance of the most ambitious improvement goal for \textbf{RTI-mean} and \textbf{RT-CVaR} (Figures~\ref{TestFigm3} and~\ref{TestFigcv3}) is highly variable between the regions. In particular, we find that some regions (e.g., Toronto, Simcoe, Durham, Hamilton) perform worse on the $\gamma=3$ testing sets than the corresponding training set. On the other hand, the test set performance for $\gamma=1$ and $\gamma=2$ is comparable to the training set performance. We believe that the observed degradation in test set performance for $\gamma=3$ is due to the model overfitting the training set for the most ambitious improvement goals. Previous research has shown that OHCA locations exhibit temporal stability \citep{Chan2016}, meaning that the locations do not vary significantly over time. Since our training and test set OHCAs are generated using the same models, the exact OHCA locations mimic the temporal stability found in practice. However, for certain regions, practitioners may prefer to use the $\gamma=2$ or 30\% improvement solutions instead of $\gamma=3$ or 50\% solutions because the model appears less likely to overfit.

\begin{figure}[t]
\begin{center}
\subfigure[\ Toronto -- \textbf{RTI-mean} \label{DistFigmToronto}]{
\frame{\includegraphics[width=.32\textwidth]{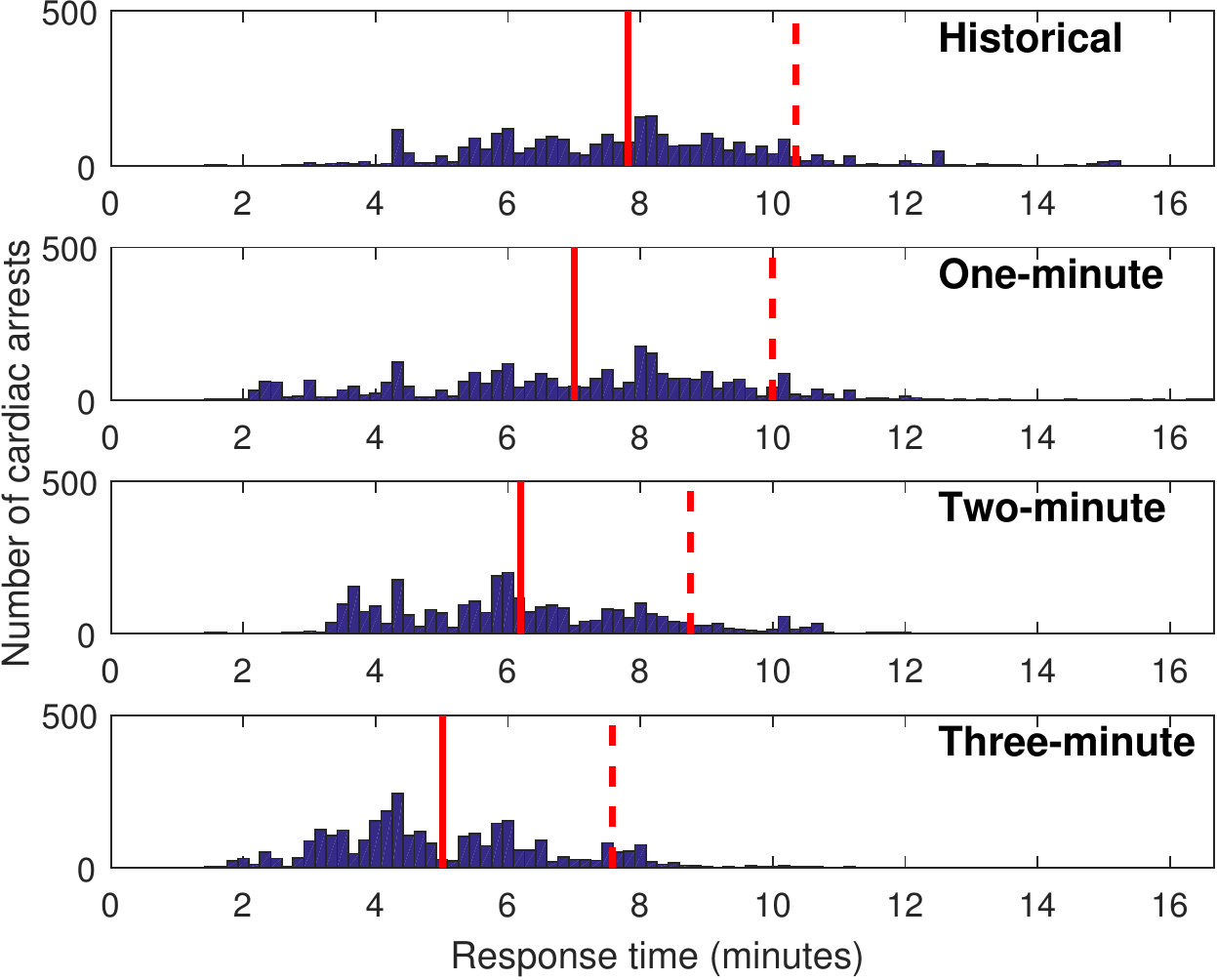}}}
\subfigure[\ Peel -- \textbf{RTI-mean} \label{DistFigmPeel}]{
\frame{\includegraphics[width=.32\textwidth]{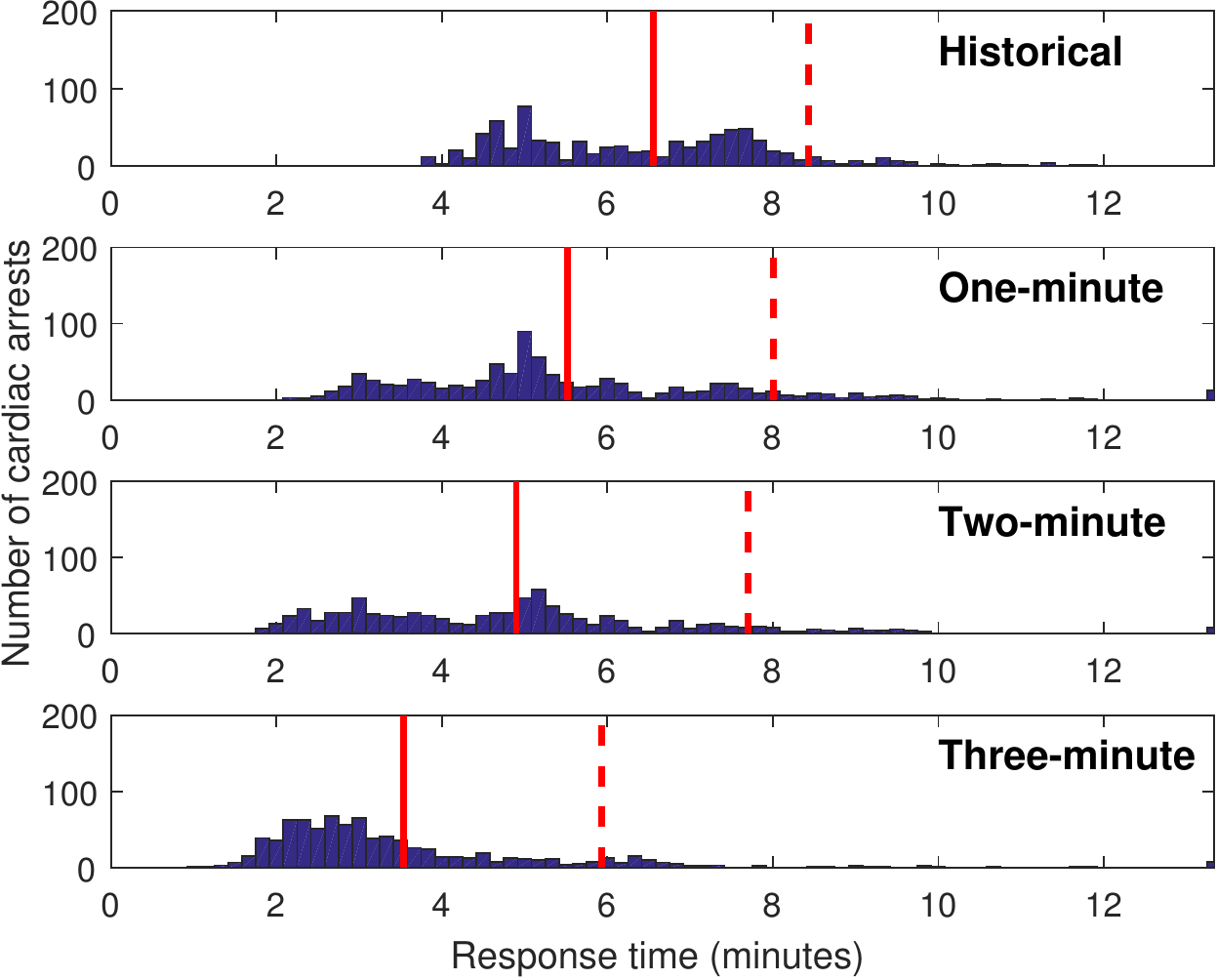}}}
\subfigure[\ Muskoka -- \textbf{RTI-mean} \label{DistFigmMuskoka}]{
\frame{\includegraphics[width=.32\textwidth]{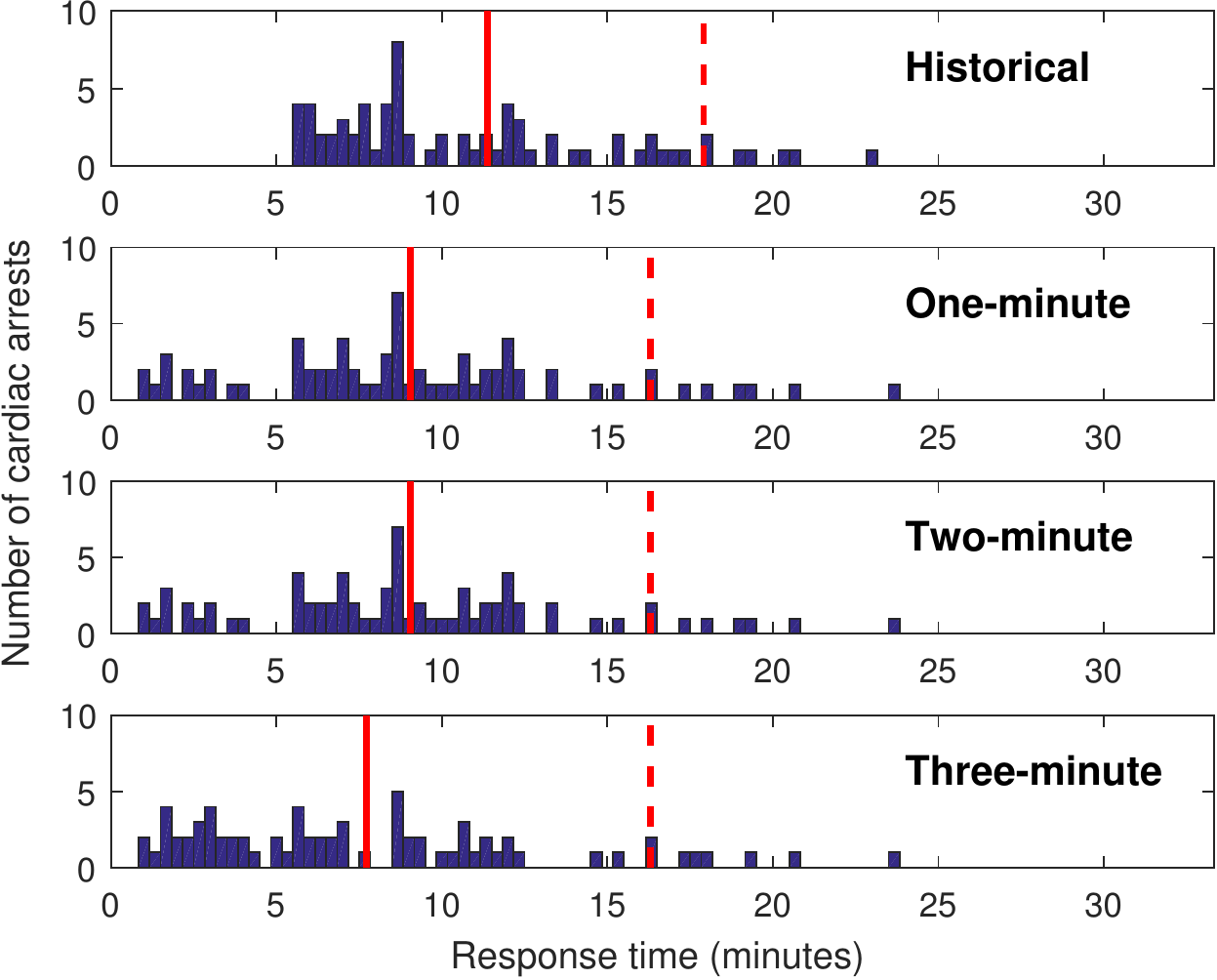}}}
\subfigure[\ Toronto -- \textbf{RT-CVaR} \label{DistFigcvToronto}]{
\frame{\includegraphics[width=.32\textwidth]{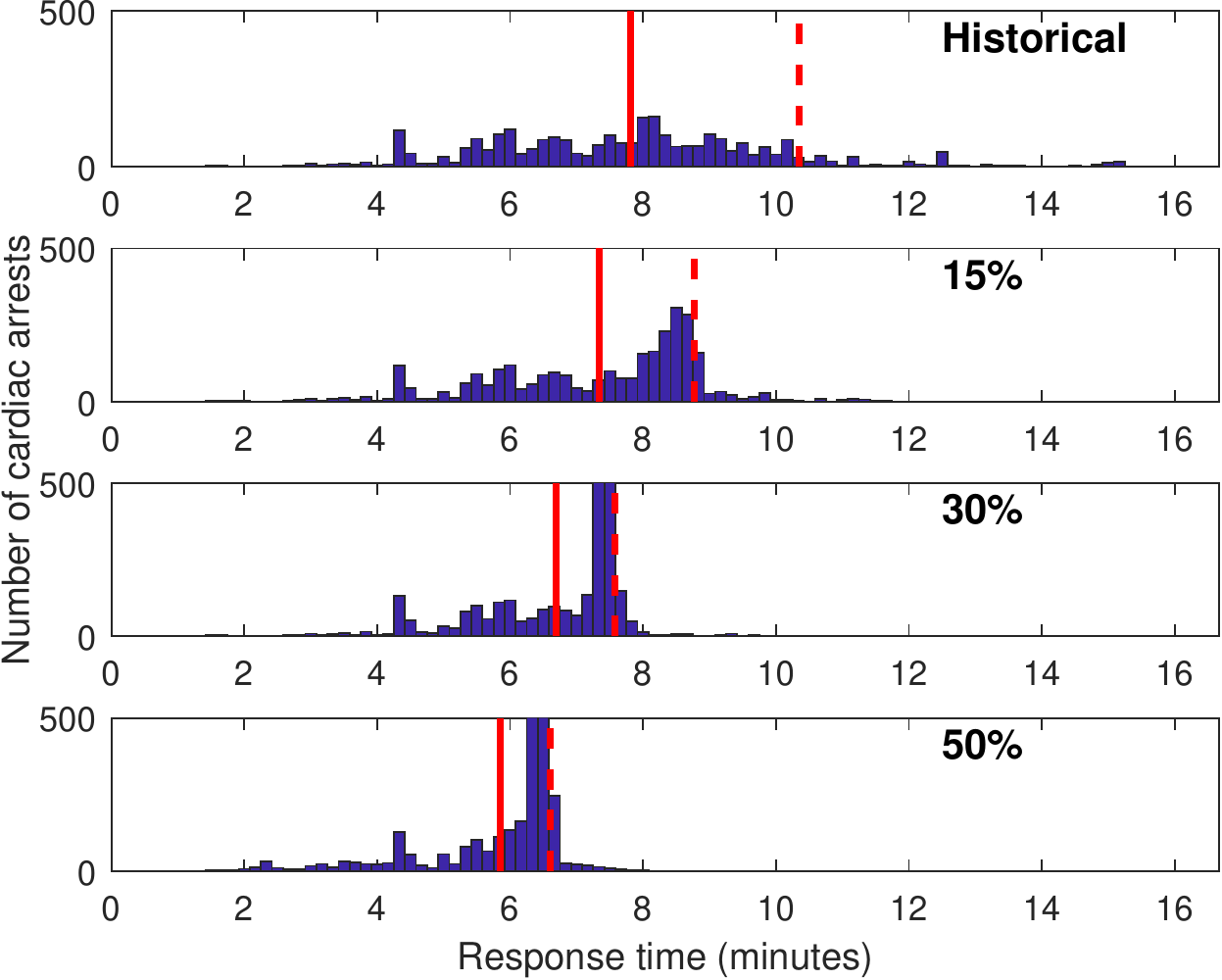}}}
\subfigure[\ Peel -- \textbf{RT-CVaR} \label{DistFigcvPeel}]{
\frame{\includegraphics[width=.32\textwidth]{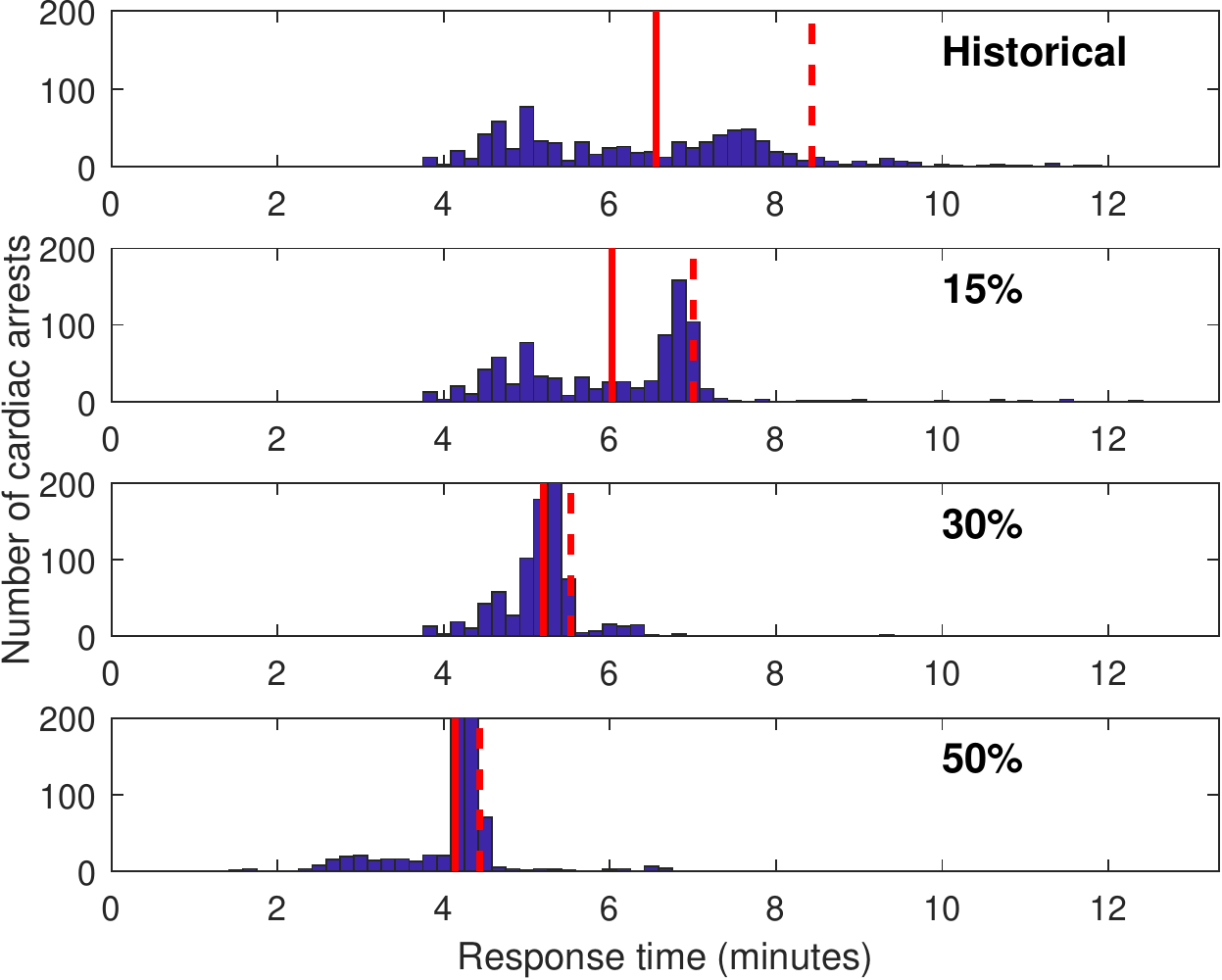}}}
\subfigure[\ Muskoka -- \textbf{RT-CVaR}\label{DistFigcvMuskoka}]{
\frame{\includegraphics[width=.32\textwidth]{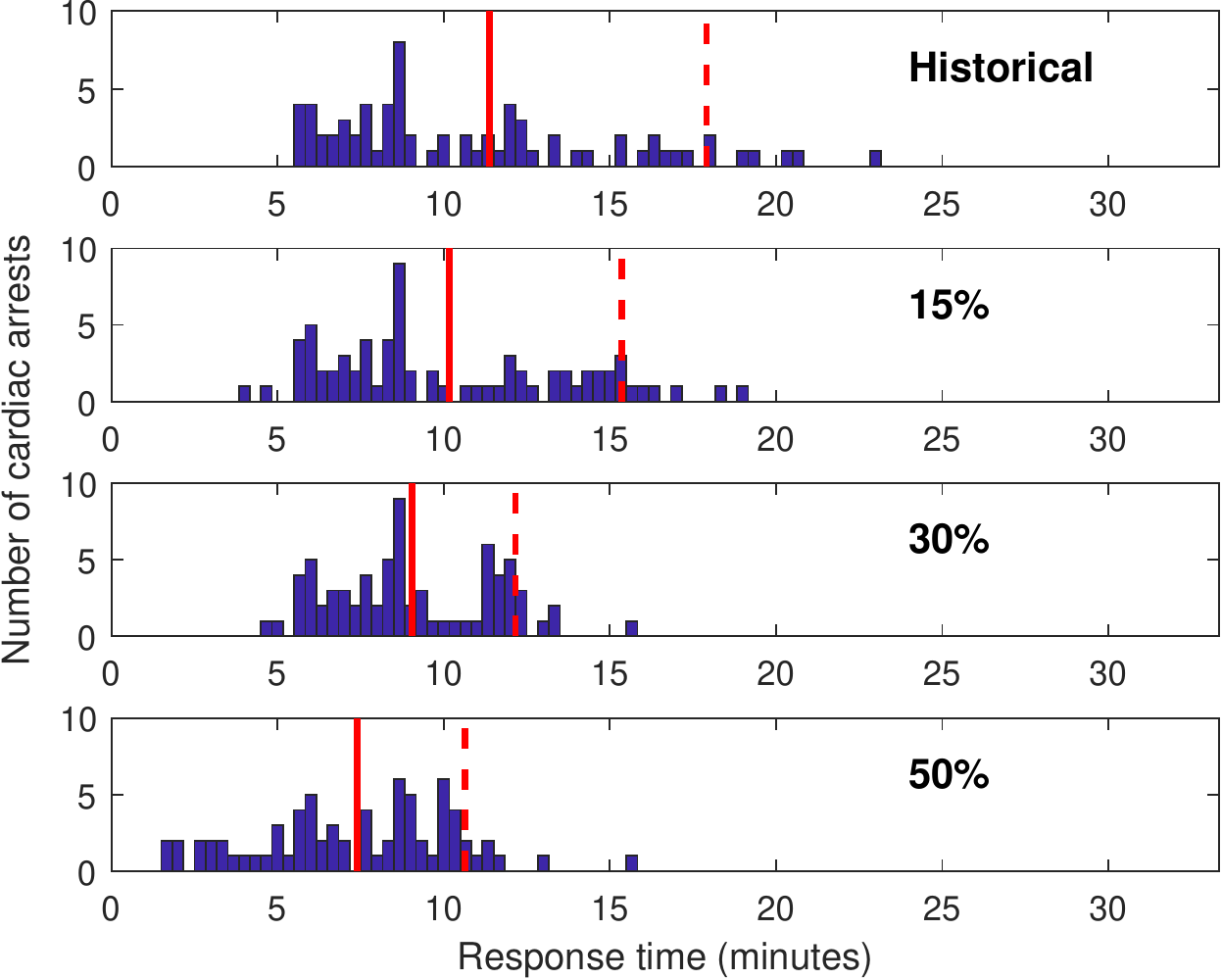}}}
\caption{Comparison of response time distributions for a random test set instance. \label{Distfig}}
\end{center}
\end{figure}

Figure~\ref{Distfig} compares the historical 911 response time distribution to the estimated response time distribution of \textbf{RTI-mean} and \textbf{RT-CVaR} in Toronto, Peel, and Muskoka. There are key differences between the response time distributions of \textbf{RTI-mean} and \textbf{RT-CVaR} that suggest that \textbf{RT-CVaR} should be used in practice. First, although \textbf{RTI-mean} is able to reduce the average response time, there is little impact on the tail of the distribution. In contrast, \textbf{RT-CVaR} is able to significantly reduce both the average and tail of the distribution, effectively reducing the variability in response times. Second, \textbf{RT-CVaR} achieves similar reductions in the average response time as compared to \textbf{RTI-mean}. Third, for the 50\% reduction in the 90th percentile of Toronto and Peel (the bottom panel in Figures~\ref{DistFigcvToronto}, and~\ref{DistFigcvPeel}), the optimized drone network is able to shrink the gap between the average and 90th percentile to within 1 minute, while shifting the entire distribution to the left. These improvements are not without cost, as \textbf{RT-CVaR} results is significantly larger drone networks and increased computational complexity to solve. 

Overall, these results provide three practical managerial insights. First, optimizing for tail performance (e.g., through \textbf{RT-CVaR}) appears to be better suited for use in practice because such an approach also significantly improves average performance, while providing more equitable coverage by locating drone bases in both cities and rural areas. Second, overfitting is a real concern when designing drone networks for ambitious improvement goals. Even though OHCAs exhibit temporal stability, practitioners may prefer to focus on modest improvement goals to avoid efficiency loss due to overfitting. Third, optimizing all resources of the drone network simultaneously (bases and drones) can reduce resource needs by over 50\% without degrading performance, compared to optimizing bases and drones separately \citep{Boutilier2017}.


\subsection{Sensitivity}\label{Sens}

In this section, we conduct a sensitivity analysis on the size of the drone network as a function of the arrival rate parameter $\lambda$ (i.e., call volume) in each region. The results in Section~\ref{Network} consider a suspected call volume multiplier of five, meaning that the number of suspected OHCAs is five times the number of confirmed OHCAs. This section considers the impact of increasing and decreasing the call volume multiplier.

Table~\ref{TableSens} provides the number of drone bases and the number of total drones required for each region. Toronto sees the largest variation in the number of required drone resources as a function of call volume, while Muskoka exhibits almost no variation, with only a single extra drone required when increasing the suspected call volume multiplier from two to ten. For most instances, we find that the number of drone bases is not impacted by changes in the call volume. This result aligns with the fact that OHCA locations exhibit temporal stability. In other words, increasing the call volume produces more OHCAs, but the locations of these OHCAs are not geographically different. As a result, the model focuses on adding more drones to existing bases.

We also find that the overall size of the drone network (bases and total drones) is insensitive to call volume for the large rural regions (e.g., Muskoka, Simcoe, York), especially for \textbf{RT-CVaR}. This insensitivity occurs because these regions require many drone bases spread over a large area to meet the response time constraints, and these drone bases are not fully utilized due to the low number of OHCAs. As a result, increasing or decreasing the call volume does not impact the number of drone resources, but rather, modifies the utilization level at each base.

\renewcommand{\arraystretch}{1.2}
\begin{table}[h!]
\caption{Sensitivity analysis on drone network size (Bases $\mid$ Total drones) for \textbf{RTI-mean} and \textbf{RT-CVaR}.}\label{TableSens}
\centering
\resizebox{\textwidth}{!}{%
\begin{tabular}{c c | c c c c c c c c}
\noalign{\smallskip} \hline
Suspected call && \multicolumn{8}{c}{\textbf{Region}} \\ \cline{3-10}
volume multiplier & Improvement & Toronto & Peel & Simcoe & York & Halton & Durham & Hamilton & Muskoka  \\ \hline
&&\multicolumn{8}{c}{\textbf{RTI-mean}} \\ \hline
\multirow{3}{*}{$2$}
&\emph{$\gamma=1$} & 1 $\mid$ 3 & 1 $\mid$ 2 & 1 $\mid$ 2 & 1 $\mid$ 2 & 2 $\mid$ 2 & 1 $\mid$ 2 & 1 $\mid$ 2 & 1 $\mid$ 1 \\
&\emph{$\gamma=2$} & 2 $\mid$ 5 & 2 $\mid$ 4 & 3 $\mid$ 4 & 2 $\mid$ 4 & 3 $\mid$ 3 & 2 $\mid$ 3 & 1 $\mid$ 2 & 1 $\mid$ 1\\
&\emph{$\gamma=3$} & 3 $\mid$ 8 & 5 $\mid$ 8 & 5 $\mid$ 6 & 5 $\mid$ 6 & 4 $\mid$ 5 & 4 $\mid$ 5 & 2 $\mid$ 3 & 2 $\mid$ 2 \\
 \hline
\multirow{3}{*}{$10$}
&\emph{$\gamma=1$} & 2 $\mid$ 6 & 1 $\mid$ 3 & 1 $\mid$ 2 & 1 $\mid$ 3 & 1 $\mid$ 2 & 1 $\mid$ 3 & 1 $\mid$ 3 & 1 $\mid$ 1\\
&\emph{$\gamma=2$} & 4 $\mid$ 12 & 3 $\mid$ 7 & 4 $\mid$ 6 & 2 $\mid$ 5 & 2 $\mid$ 5 & 2 $\mid$ 5 & 2 $\mid$ 4 & 2 $\mid$ 2\\
&\emph{$\gamma=3$} & 6 $\mid$ 18 & 5 $\mid$ 12 & 6 $\mid$ 10 & 3 $\mid$ 9 & 4 $\mid$ 8 & 3 $\mid$ 8 & 2 $\mid$ 6 & 2 $\mid$ 3\\
 \hline
&&\multicolumn{8}{c}{\textbf{RT-CVaR}} \\ \hline
\multirow{3}{*}{$2$}
&15\% & 2 $\mid$ 3 & 4 $\mid$ 4 & 4 $\mid$ 4 & 4 $\mid$ 4 & 4 $\mid$ 4 & 7 $\mid$ 7 & 2 $\mid$ 2 & 2 $\mid$ 2\\
&30\% & 3 $\mid$ 5 & 7 $\mid$ 7 & 6 $\mid$ 6 & 6 $\mid$ 6 & 5 $\mid$ 5 & 12 $\mid$ 12 & 3 $\mid$ 3 & 4 $\mid$ 4\\
&50\% & 8 $\mid$ 13 & 13 $\mid$ 13 & 17 $\mid$ 17 & 13 $\mid$ 13 & 9 $\mid$ 9 & Infeasible & 6 $\mid$ 7 & 6 $\mid$ 6 \\
\hline
\multirow{3}{*}{$10$}
&15\% & 2 $\mid$ 5 & 5 $\mid$ 6 & 4 $\mid$ 4 & 5 $\mid$ 5 & 4 $\mid$ 4 & 8 $\mid$ 10 & 2 $\mid$ 3 & 2 $\mid$ 2 \\
&30\% & 3 $\mid$ 9 & 9 $\mid$ 11 & 7 $\mid$ 7 & 7 $\mid$ 9 & 6 $\mid$ 7 & 13 $\mid$ 16 & 3 $\mid$ 5 & 4 $\mid$ 4\\
&50\% & 9 $\mid$ 23 & 14 $\mid$ 20 & 18 $\mid$ 18 & 14 $\mid$ 18 & 13 $\mid$ 15 & Infeasible & 6 $\mid$ 10 & 6 $\mid$ 6\\
\hline
\end{tabular}}
\end{table}

\subsection{The trade-off between bases and drones}\label{DBtrade}

\begin{figure}[t]
\begin{center}
\subfigure[\ Toronto -- \textbf{RTI-mean}  \label{AlphaM1}]{
\frame{\includegraphics[width=.3\textwidth]{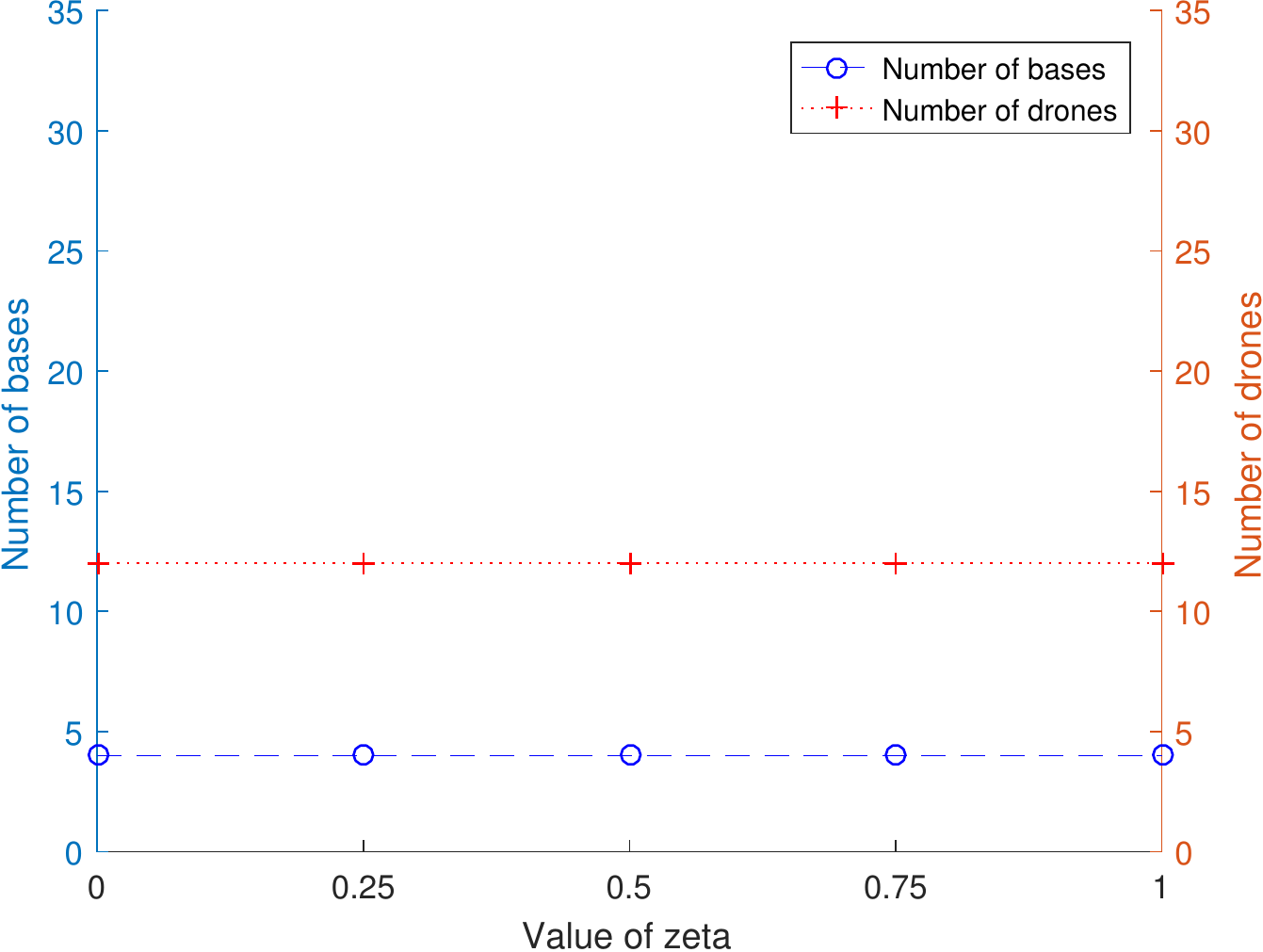}}}
\subfigure[\ Peel -- \textbf{RTI-mean}  \label{AlphaM5}]{
\frame{\includegraphics[width=.3\textwidth]{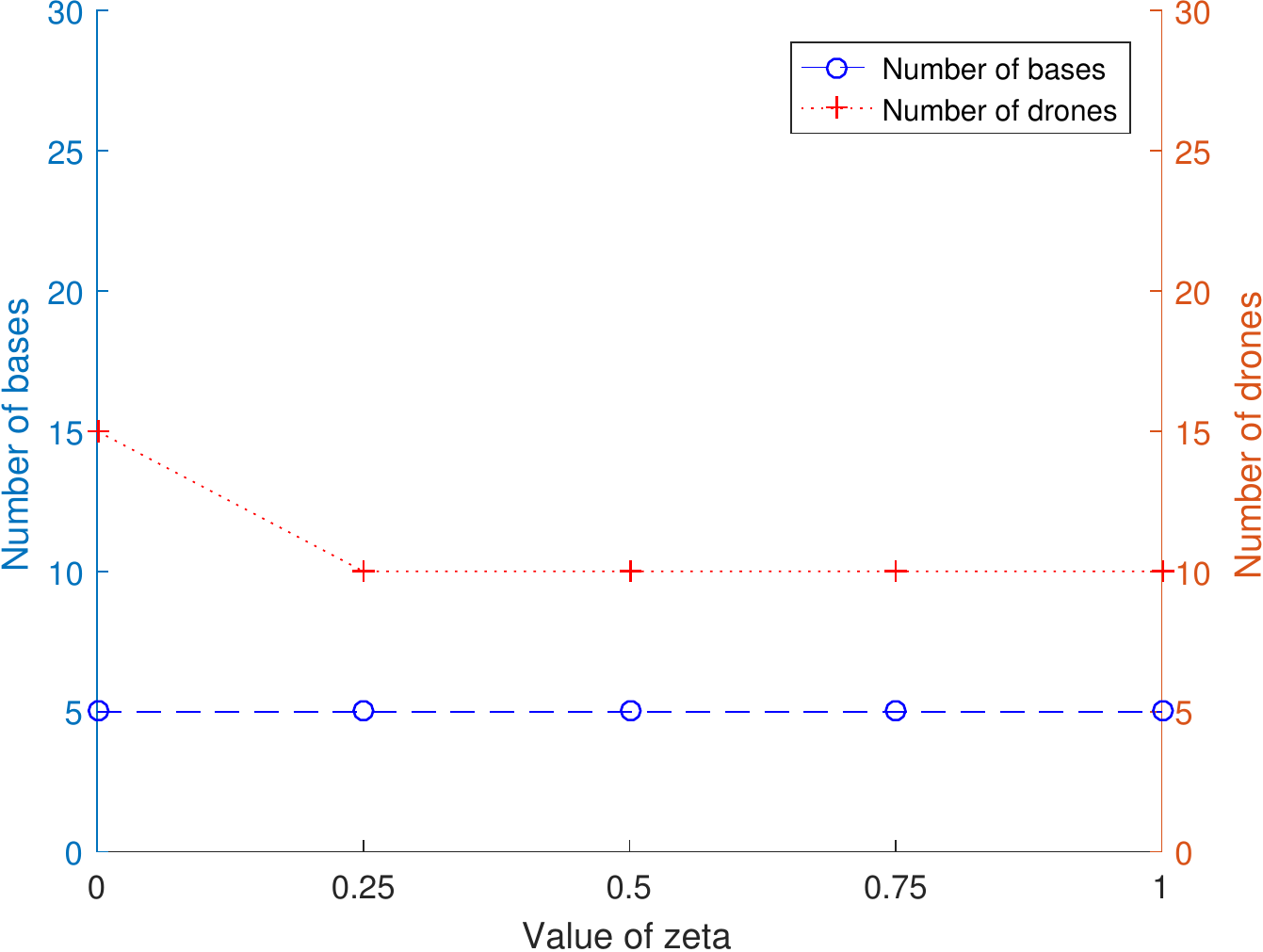}}}
\subfigure[\ Muskoka -- \textbf{RTI-mean}  \label{AlphaM4}]{
\frame{\includegraphics[width=.3\textwidth]{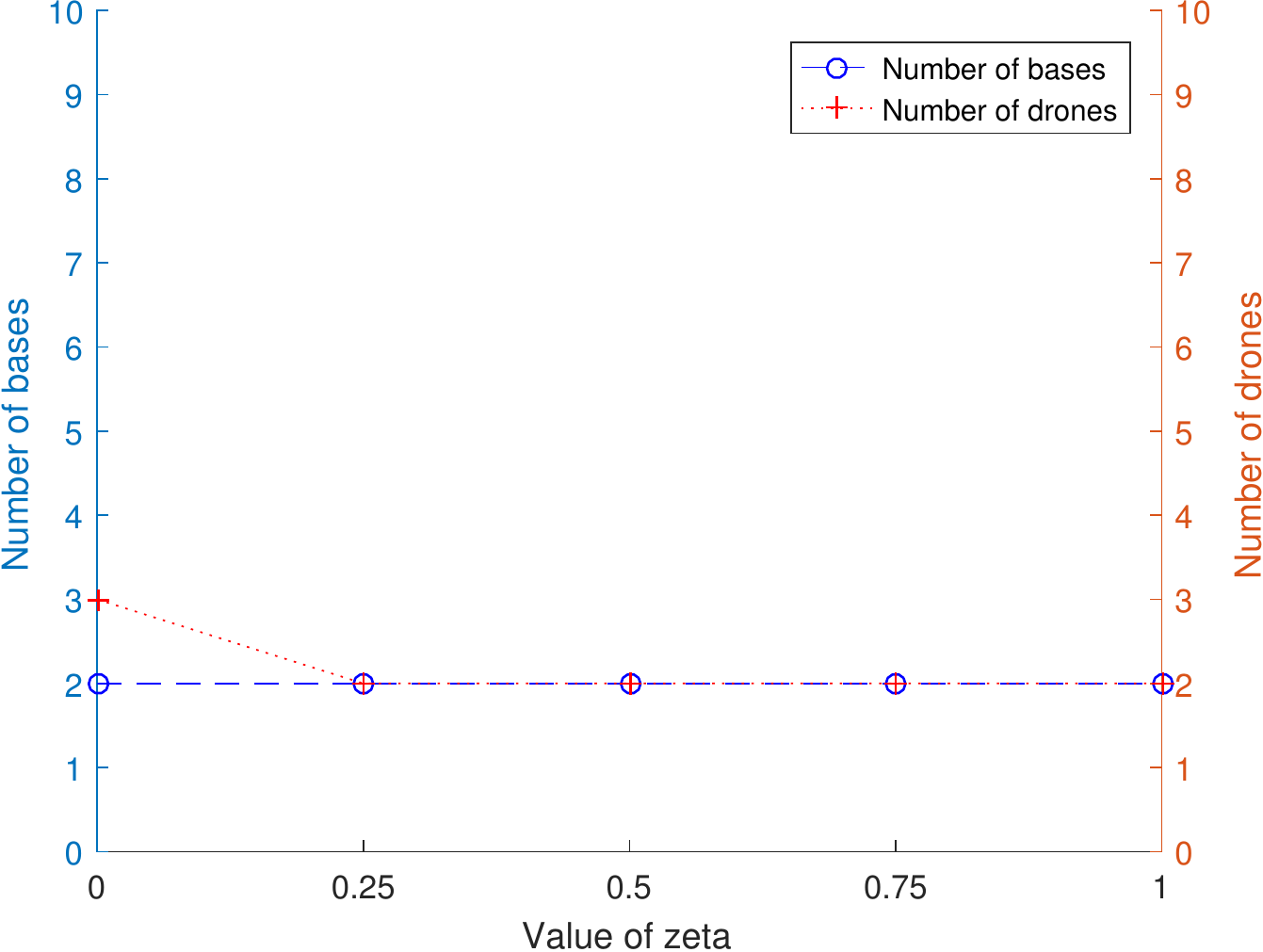}}}
\subfigure[\ Toronto -- \textbf{RT-CVaR} \label{AlphaC1}]{
\frame{\includegraphics[width=.3\textwidth]{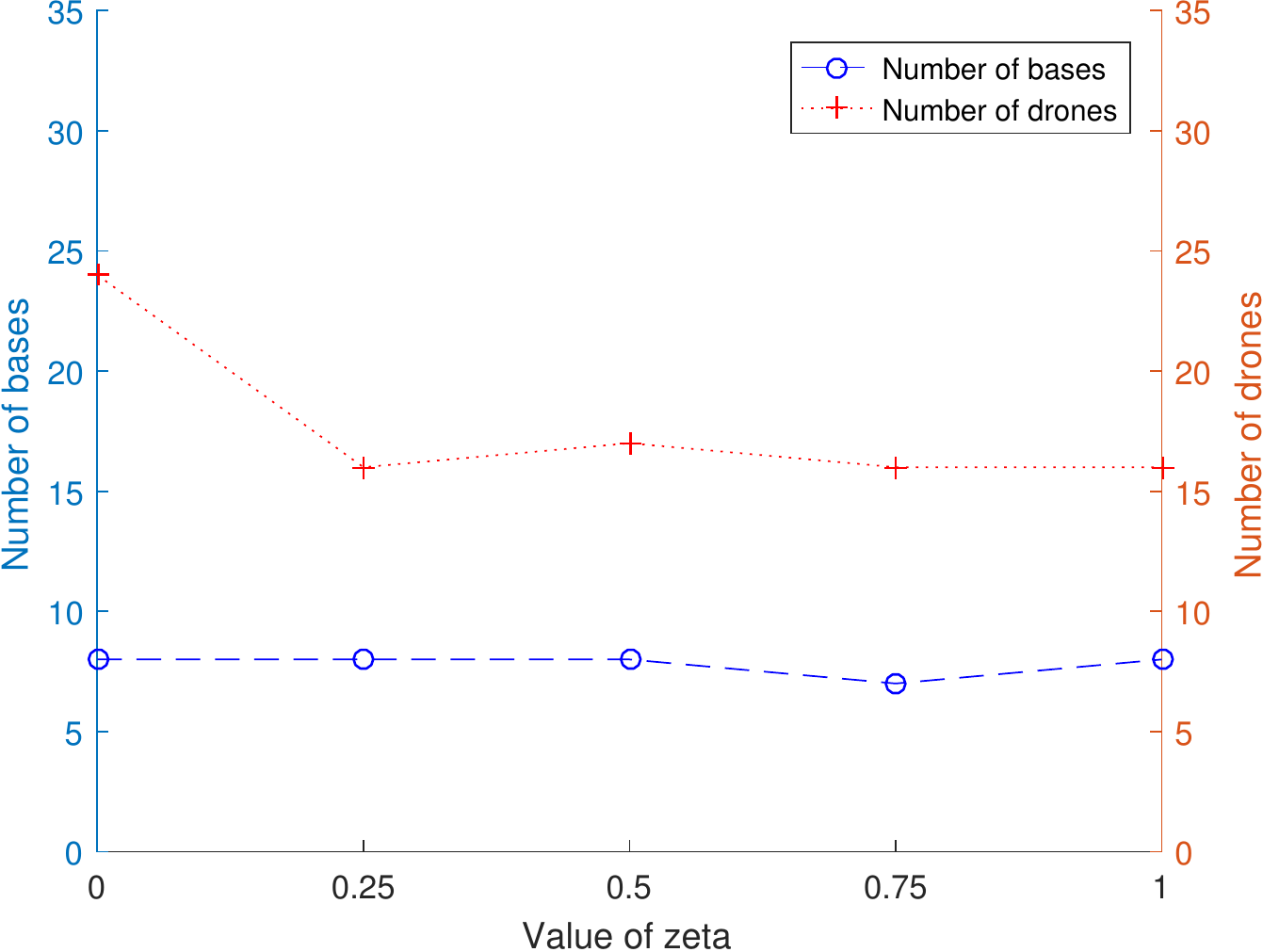}}}
\subfigure[\ Peel -- \textbf{RT-CVaR} \label{AlphaC5}]{
\frame{\includegraphics[width=.3\textwidth]{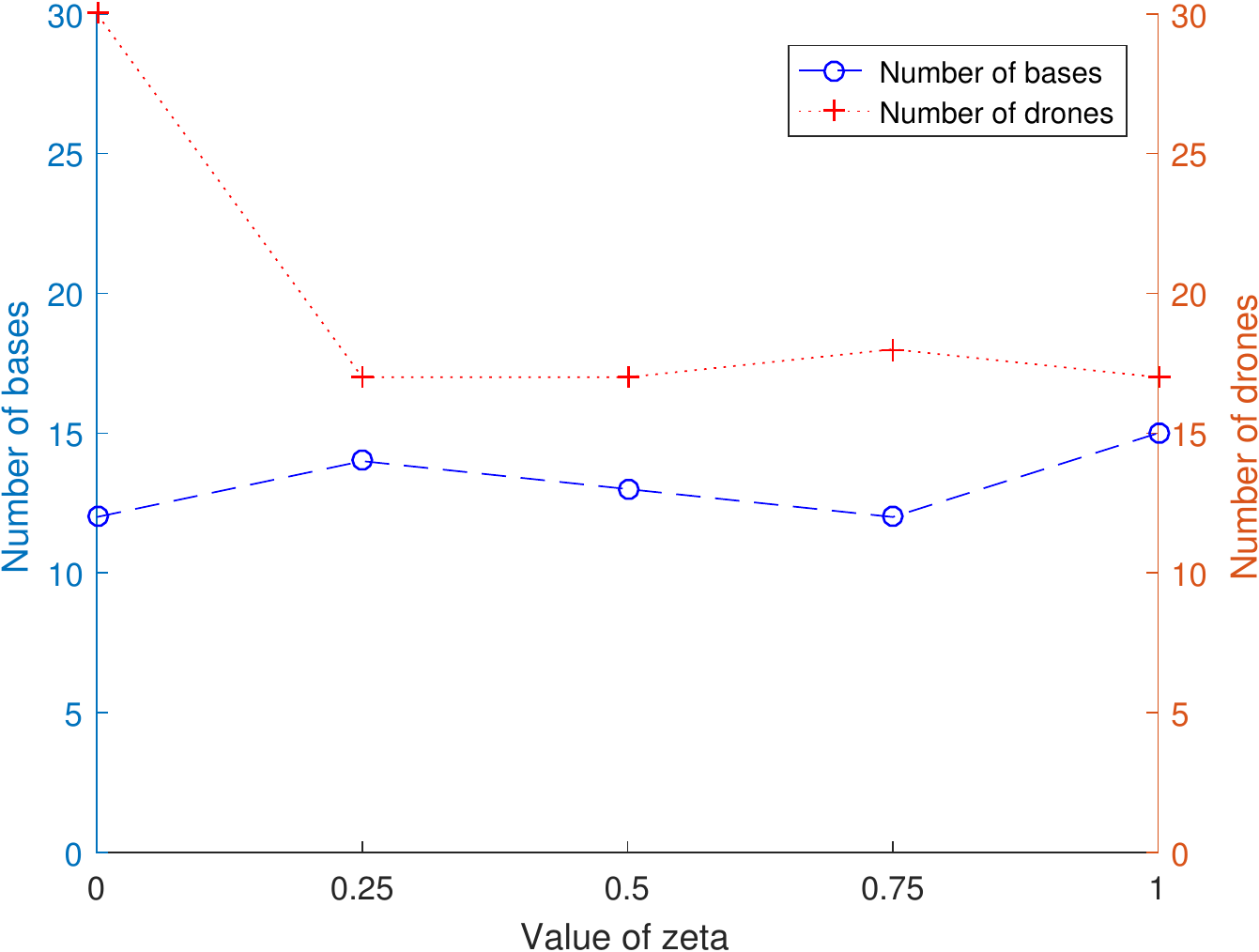}}}
\subfigure[\ Muskoka -- \textbf{RT-CVaR}\label{AlphaC4}]{
\frame{\includegraphics[width=.3\textwidth]{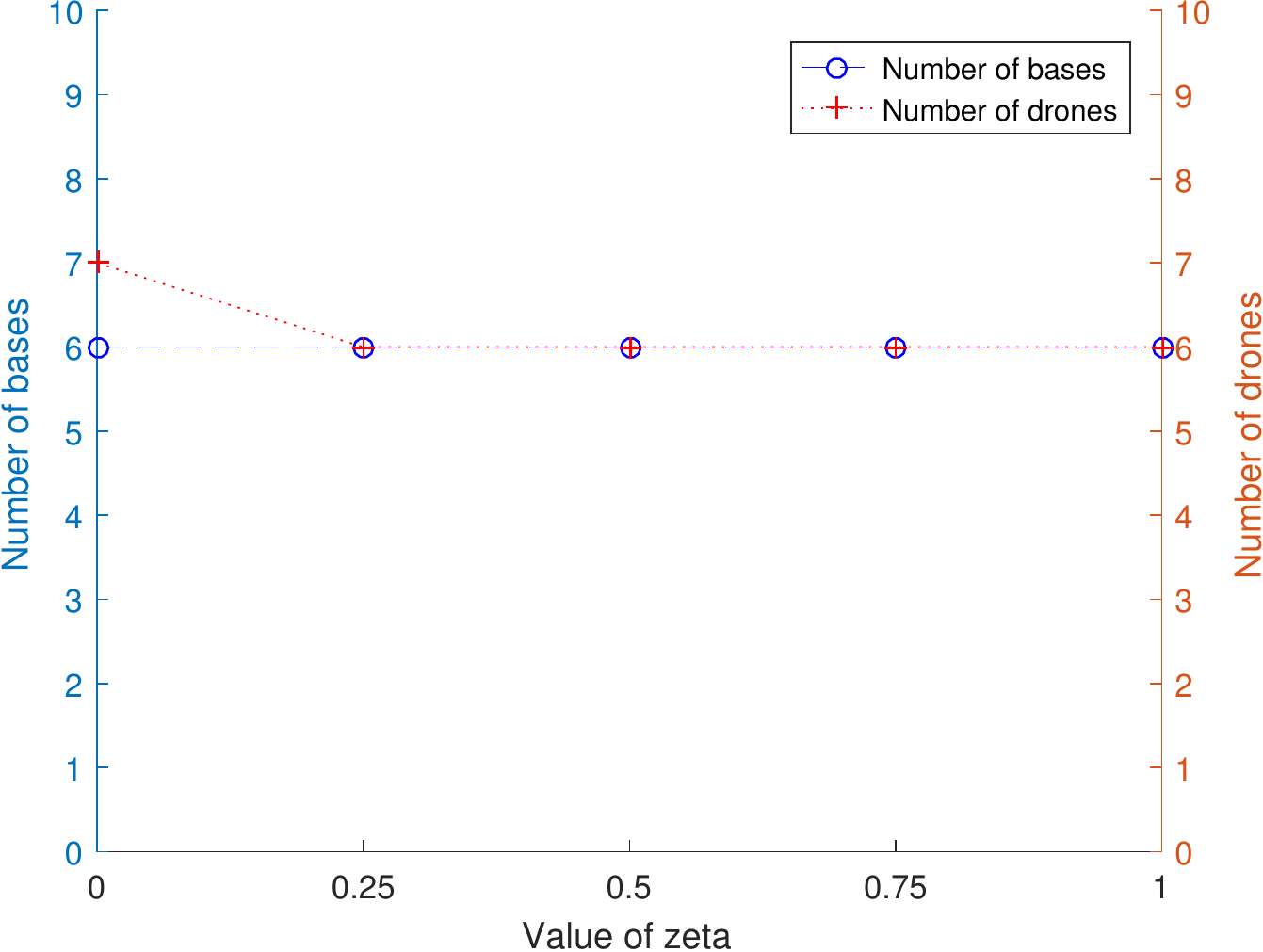}}}
\caption{Comparison of the trade-off between bases and drones.\label{AlphaFig}}
\end{center}
\end{figure}

In this section, we present results that investigate the trade-off between minimizing the total number of drones versus minimizing the number of drones bases. Figure~\ref{AlphaFig} displays the number of drone bases and the total number of drones for \textbf{RTI-mean} and \textbf{RT-CVaR} in Toronto, Peel, and Muskoka. With the exception of $\zeta=0$, there is little impact on the total number of drones or the number of bases, suggesting that our results are not influenced by modifying the objective to focus on the number of total drones or bases. When $\zeta=0$, the objective focuses only on minimizing the number of bases, with no penalty for additional drones, leading to solutions with an artificially large number of drones. The practical implication of this result is important: the optimal drone network is insensitive to cost differences between building new bases or adding additional drones to current bases. As a result, EMS providers may not need to consider the cost trade-off associated with centralizing drone resources, regardless of geographic region.

\section{Discussion}\label{Disc}

In this section, we summarize our work (Section~\ref{Sum}), present the results through an urban vs. rural lens (Section~\ref{UR}), highlight the potential benefits of using drones to deliver AEDs (Section~\ref{Benefits}), and discuss the key factors associated with implementation (Section~\ref{Implement}).

\subsection{Summary}\label{Sum}

In this paper, we developed an integrated location-queuing model that is based on the $p$-median architecture, where each base constitutes an explicit $M/M/d$ queue, and that incorporates estimated baseline response times to the demand points. We further generalized our model by formulating a conditional value-at-risk variant that allows us to control the tail of the response time distribution. We then developed a reformulation technique that exploits the baseline response times to induce coverage-like sparsity for median-type problems, allowing us to solve large-scale, real-world instances to optimality using an off-the-shelf solver. Both our model and reformulation technique are generalizable to any location-queuing problem where transportation options are added on top of an existing system to optimize performance.
 
We demonstrated the application of our framework to determine the optimal deployment of AED-enabled drone resources to meet various response time targets using eight years of real data from an area covering 26,000 square kilometers around Toronto, Canada. In the process, we developed a two-stage machine learning approach to simulate cardiac arrest incidents. Our experiments indicate that a modest number of drone resources are able to significantly reduce response times. Furthermore, we demonstrated that for most regions (except the most urban), the number of drone bases is insensitive to the call volume and for all regions, the resulting drone network does not depend on the cost difference between building more bases or adding more drones to existing bases. Overall, drone-delivered AEDs have the potential to be a transformative innovation in the provision of emergency care and this paper provides a realistic framework that can be leveraged by system designers and/or EMS personnel seeking to investigate design questions associated with a drone network.

\subsection{Viewing the results through an urban vs. rural lens}\label{UR}

As outlined in the Introduction, rural areas have slower response times and lower survival rates from OHCA as compared to urban areas. Our results demonstrate that drones have the ability to significantly improve the response time distribution for rural areas, effectively eliminating the response time difference between urban and rural areas. However, rural areas also have fewer OHCAs, meaning that drone utilization will be limited as compared to urban areas. In other words, a drone network in a rural area will have a large impact on a small number of OHCAs, while a drone network in an urban area will have a small impact on a large number of OHCAs. We summarize other key differences between rural and urban drone networks: 
\begin{itemize}
\item Optimizing the average response time prescribed solutions that locate drone bases in urban areas, while optimizing the 90th percentile prescribed solutions that locate drone bases in both urban and rural areas.
\item Rural regions require more drone bases, but fewer drones per base compared to urban regions.
\item Rural regions exhibit a higher variability in test set performance most likely due to fewer OHCAs, especially when optimizing for the tail of the response time distribution.
\item Rural regions are less sensitive to changes in call volume as compared to urban regions, suggesting greater utilization of drone resources in urban areas.
\end{itemize}
Overall, these observations suggest that an equity-efficiency trade-off between drone networks in rural and urban areas will be a key consideration for system designers implementing a drone network in conjunction with EMS.

\subsection{Potential benefits}\label{Benefits}

There are many potential benefits of using a drone network to deliver AEDs. First and foremost, drone-delivered AEDs have the potential to improve survival for patients with OHCA because the probability of survival is a decreasing function of response time \citep{Valenzuela1997}.

Second, a drone network has the potential to improve AED utilization, which is critical for improving survival. Currently, an AED is used in less than 3\% of all OHCAs \citep{Weisfeldt2010, Boutilier2017}, which can be partially attributed to access and availability issues \citep{Sun2016}, especially in rural areas. Drone-delivered AEDs effectively eliminate access and availability issues by delivering an AED to the scene of all OHCAs, regardless of location or time of day. In contrast to public access defibrillation programs, which deploy static AEDs in the community, drones actively mobilize AEDs for both public and private location OHCAs. Moreover, the drone is equipped with a microphone and camera that can be leveraged by the 911 dispatcher to survey the scene and provide assistance to bystanders, which has been shown to increase participation and AED utilization \citep{Lerner2012}. 

Third, drones may be able to deliver AEDs to OHCAs in hard-to-reach areas such as golf courses, hiking trails, and balconies of high-rise buildings, which are known to suffer a survival disadvantage \citep{Drennan2015}. Drones are able to provide rapid response to these areas due to straight-line flight, which avoids traffic and does not require a road network.

\subsection{Implementation factors}\label{Implement}

Before a drone network can be realized and implemented, there are several operational, regulatory, technical, and educational challenges that must be addressed. Most importantly, EMS providers and policy makers must determine what type of drone network is desired. Practical questions such as whether an EMS system should focus on improving the average or 90th percentile in either rural or urban areas can be investigated using our models.
	
Operational decisions will influence how the resulting network is implemented. The most critical question pertains to how the drone system will be integrated with the current EMS system. EMS providers will need to decide which calls drones will be dispatched to, how the drone will be transported back to base after use, and how to conduct regular and after-use maintenance of the drones. Our sensitivity analysis explores the impact of call volume on the size of the drone network, which is impacted by each of these operational considerations. Our main findings suggest that the number and location of drone bases are not impacted by call volume, but the number of drones allocated to each base is.

Many regulatory challenges must be addressed before drone-delivered AEDs can be realized. Currently, the most stifling regulation limits drone use to operator line of sight, meaning that they must be flown where a human operator can see them. However, significant lobbying efforts by major corporations are paving the way for more progressive legislation. For example, countries like Canada and the USA have begun handling drone regulations on a case-by-case basis. In Canada, approval has been granted to conduct a pilot study to test the efficacy of drone-delivered AEDs. Pilot studies can leverage our models to determine the most effective location(s) to test drone-delivered AEDs.

There are also technical hurdles that require further advancement to promote the safe operations of drone-delivered AEDs. Much research has focused on improving computer vision, autonomous flight, and object avoidance measures, which are critical for drone deliveries. Continuous improvement in these areas combined with advancements in bad weather flight, battery capacity, and payload size will pave the way for future EMS drone applications. 

Aside from operational, technical, and regulatory issues, a drone network must also be accompanied by widespread educational campaigns. In particular, educational campaigns should focus on teaching the public about the signs and symptoms of OHCA, and the importance of initializing EMS response. Further education is needed to alert the public that EMS drones are a life-saving device and should not be tampered with. Similar to ambulances, uniformly coloured drones with clear EMS markings, lights, and sirens are imperative to achieving this goal. A tangential benefit of drone-delivered AEDs is the opportunity to create public excitement and enthusiasm, which can further aid with education and awareness around OHCA. 



%
%
%


\bibliographystyle{informs2014} 
\bibliography{testbib} 

\begin{thebibliography}{62}
\providecommand{\natexlab}[1]{#1}
\providecommand{\url}[1]{\texttt{#1}}
\providecommand{\urlprefix}{URL }

\bibitem[{Aboolian et~al.(2008)Aboolian, Berman, \protect\BIBand{}
  Drezner}]{Aboolian2008}
Aboolian R, Berman O, Drezner Z (2008) Location and allocation of service units
  on a congested network. \emph{IIE Transactions} 40(4):422--433,
  \urlprefix\url{http://dx.doi.org/10.1080/07408170701411385}.

\bibitem[{Agatz et~al.(2018)Agatz, Bouman, \protect\BIBand{}
  Schmidt}]{Agatz2018}
Agatz N, Bouman P, Schmidt M (2018) Optimization approaches for the traveling
  salesman problem with drone. \emph{Transportation Science}
  \urlprefix\url{http://dx.doi.org/10.1287/trsc.2017.0791}.

\bibitem[{Ahmadi-Javid et~al.(2017)Ahmadi-Javid, Seyedi, \protect\BIBand{}
  Syam}]{Ahmadi2017}
Ahmadi-Javid A, Seyedi P, Syam SS (2017) A survey of healthcare facility
  location. \emph{Computers and Operations Research} 79:223--263.

\bibitem[{Altman(1992)}]{Altman1992}
Altman NS (1992) An introduction to kernel and nearest-neighbor nonparametric
  regression. \emph{The American Statistician} 46(3):175--185,
  \urlprefix\url{http://dx.doi.org/10.1080/00031305.1992.10475879}.

\bibitem[{Basar et~al.(2012)Basar, Catay, \protect\BIBand{}
  Unluyurt}]{Basar2012}
Basar A, Catay B, Unluyurt T (2012) A taxonomy for emergency service station
  location problem. \emph{Optimization letters} 6:1147--1160.

\bibitem[{Batta \protect\BIBand{} Berman(1989)}]{Batta1989}
Batta R, Berman O (1989) A location model for a facility operating as an m/g/k
  queue. \emph{Networks} 19(6):717--728,
  \urlprefix\url{http://dx.doi.org/10.1002/net.3230190609}.

\bibitem[{Berman \protect\BIBand{} Drezner(2007)}]{Berman2007}
Berman O, Drezner Z (2007) The multiple server location problem. \emph{Journal
  of the Operational Research Society} 58(1):91--99, ISSN 1476-9360,
  \urlprefix\url{http://dx.doi.org/10.1057/palgrave.jors.2602126}.

\bibitem[{Berman et~al.(1985)Berman, Larson, \protect\BIBand{}
  Chiu}]{Berman1985}
Berman O, Larson RC, Chiu SS (1985) Optimal server location on a network
  operating as an m/g/1 queue. \emph{Operations Research} 33(4):746--771,
  \urlprefix\url{http://dx.doi.org/10.1287/opre.33.4.746}.

\bibitem[{Bishop(2006)}]{Bishop2006}
Bishop C (2006) \emph{Pattern Recognition and Machine Learning} (Sprnger).

\bibitem[{Boutilier et~al.(2017)Boutilier, Brooks, Janmohamed, Byers, Buick,
  Zhan, Schoellig, Cheskes, Morrison, \protect\BIBand{} Chan}]{Boutilier2017}
Boutilier JJ, Brooks SC, Janmohamed A, Byers A, Buick JE, Zhan C, Schoellig AP,
  Cheskes S, Morrison LJ, Chan TCY (2017) Optimizing a drone network to deliver
  automated external defibrillatorsclinical perspective. \emph{Circulation}
  135(25):2454--2465, ISSN 0009-7322,
  \urlprefix\url{http://dx.doi.org/10.1161/CIRCULATIONAHA.116.026318}.

\bibitem[{Brotcorne et~al.(2003)Brotcorne, Laporte, \protect\BIBand{}
  Semet}]{Brotcorne2003}
Brotcorne L, Laporte G, Semet F (2003) Ambulance location and relocation
  models. \emph{European Journal of Operational Research} 147(3):451 -- 463,
  ISSN 0377-2217,
  \urlprefix\url{http://dx.doi.org/https://doi.org/10.1016/S0377-2217(02)00364-8}.

\bibitem[{Caffrey et~al.(2002)Caffrey, Willoughby, Pepe, \protect\BIBand{}
  Becker}]{Caffrey2002}
Caffrey SL, Willoughby PJ, Pepe PE, Becker LB (2002) Public use of automated
  external defibrillators. \emph{New England Journal of Medicine}
  347(16):1242--1247, \urlprefix\url{http://dx.doi.org/10.1056/NEJMoa020932},
  pMID: 12393821.

\bibitem[{Canada(2018)}]{TransportCanada}
Canada T (2018) Flying your drone safely and legally. Technical report,
  Transport Canada.

\bibitem[{Carlsson \protect\BIBand{} Song(2018)}]{Carlsson2018}
Carlsson JG, Song S (2018) Coordinated logistics with a truck and a drone.
  \emph{Management Science}
  \urlprefix\url{http://dx.doi.org/10.1287/mnsc.2017.2824}.

\bibitem[{Carson \protect\BIBand{} Batta(1990)}]{Carson1990}
Carson YM, Batta R (1990) Locating an ambulance on the amherst campus of the
  state university of new york at buffalo. \emph{Interfaces} 20:43--49.

\bibitem[{Chan et~al.(2013)Chan, Li, Lebovic, Tang, Chan, Cheng, Morrison,
  \protect\BIBand{} Brooks}]{Chan2013}
Chan TC, Li H, Lebovic G, Tang SK, Chan JY, Cheng HC, Morrison LJ, Brooks SC
  (2013) Identifying locations for public access defibrillators using
  mathematical optimizationclinical perspective. \emph{Circulation}
  127(17):1801--1809, ISSN 0009-7322,
  \urlprefix\url{http://dx.doi.org/10.1161/CIRCULATIONAHA.113.001953}.

\bibitem[{Chan et~al.(2016)Chan, Demirtas, \protect\BIBand{} Kwon}]{Chan2016}
Chan TCY, Demirtas D, Kwon RH (2016) Optimizing the deployment of public access
  defibrillators. \emph{Management Science} 62(12):3617--3635,
  \urlprefix\url{http://dx.doi.org/10.1287/mnsc.2015.2312}.

\bibitem[{Chan et~al.(2018)Chan, Shen, \protect\BIBand{} Siddiq}]{Chan2018}
Chan TCY, Shen ZJM, Siddiq A (2018) Robust defibrillator deployment under
  cardiac arrest location uncertainty via row-and-column generation.
  \emph{Operations Research} 66(2):358--379,
  \urlprefix\url{http://dx.doi.org/10.1287/opre.2017.1660}.

\bibitem[{Claesson et~al.(2016)Claesson, Fredman, Svensson, Ringh, Hollenberg,
  Nordberg, Rosenqvist, Djarv, {\"O}sterberg, Lennartsson, \protect\BIBand{}
  Ban}]{Claesson2016}
Claesson A, Fredman D, Svensson L, Ringh M, Hollenberg J, Nordberg P,
  Rosenqvist M, Djarv T, {\"O}sterberg S, Lennartsson J, Ban Y (2016) Unmanned
  aerial vehicles (drones) in out-of-hospital-cardiac-arrest.
  \emph{Scandinavian Journal of Trauma, Resuscitation and Emergency Medicine}
  24(1):124, ISSN 1757-7241,
  \urlprefix\url{http://dx.doi.org/10.1186/s13049-016-0313-5}.

\bibitem[{Dorling et~al.(2017)Dorling, Heinrichs, Messier, \protect\BIBand{}
  Magierowski}]{Dorling2017}
Dorling K, Heinrichs J, Messier GG, Magierowski S (2017) Vehicle routing
  problems for drone delivery. \emph{IEEE Transactions on Systems, Man, and
  Cybernetics: Systems} 47(1):70--85, ISSN 2168-2216,
  \urlprefix\url{http://dx.doi.org/10.1109/TSMC.2016.2582745}.

\bibitem[{Dorr \protect\BIBand{} Duquette(2016)}]{FAA}
Dorr L, Duquette A (2016) Fact sheet -- small unmanned aircraft regulations.
  Technical report, Federal Aviation Administration.

\bibitem[{Drennan et~al.(2016)Drennan, Strum, Byers, Buick, Lin, Cheskes, Hu,
  \protect\BIBand{} Morrison}]{Drennan2015}
Drennan IR, Strum RP, Byers A, Buick JE, Lin S, Cheskes S, Hu S, Morrison LJ
  (2016) Out-of-hospital cardiac arrest in high-rise buildings: delays to
  patient care and effect on survival. \emph{CMAJ} 188(6):413--419, ISSN
  0820-3946, \urlprefix\url{http://dx.doi.org/10.1503/cmaj.150544}.

\bibitem[{Erkut et~al.(2007)Erkut, Ingolfsson, \protect\BIBand{} Erdo{\u
  g}an}]{Erkut2007}
Erkut E, Ingolfsson A, Erdo{\u g}an G (2007) Ambulance location for maximum
  survival. \emph{Naval Research Logistics (NRL)} 55(1):42--58,
  \urlprefix\url{http://dx.doi.org/10.1002/nav.20267}.

\bibitem[{Go et~al.(2013)Go, Mozaffarian, Roger, Benjamin, Berry, Borden,
  Bravata, Dai, Ford, Fox, Franco, Fullerton, Gillespie, Hailpern, Heit,
  Howard, Huffman, Kissela, Kittner, Lackland, Lichtman, Lisabeth, Magid,
  Marcus, Marelli, Matchar, McGuire, Mohler, Moy, Mussolino, Nichol, Paynter,
  Schreiner, Sorlie, Stein, Turan, Virani, Wong, Woo, \protect\BIBand{}
  Turner}]{Go2013}
Go AS, Mozaffarian D, Roger VL, Benjamin EJ, Berry JD, Borden WB, Bravata DM,
  Dai S, Ford ES, Fox CS, Franco S, Fullerton HJ, Gillespie C, Hailpern SM,
  Heit JA, Howard VJ, Huffman MD, Kissela BM, Kittner SJ, Lackland DT, Lichtman
  JH, Lisabeth LD, Magid D, Marcus GM, Marelli A, Matchar DB, McGuire DK,
  Mohler ER, Moy CS, Mussolino ME, Nichol G, Paynter NP, Schreiner PJ, Sorlie
  PD, Stein J, Turan TN, Virani SS, Wong ND, Woo D, Turner MB (2013) Heart
  disease and stroke statistics{\textemdash}2013 update. \emph{Circulation}
  127(1):e6--e245, ISSN 0009-7322,
  \urlprefix\url{http://dx.doi.org/10.1161/CIR.0b013e31828124ad}.

\bibitem[{H3Dynamics(2018)}]{DroneBox}
H3Dynamics (2018) Drone box. Technical report, H3Dynamics,
  \urlprefix\url{https://www.h3dynamics.com/products/drone-box/}.

\bibitem[{Hallstrom et~al.(2004)Hallstrom, Ornato, Weifeldt, Travers,
  Christenson, \protect\BIBand{} Investigators}]{PAD2004}
Hallstrom A, Ornato J, Weifeldt M, Travers A, Christenson J, Investigators
  TPADT (2004) Public-access defibrillation and survival after out-of-hospital
  cardiac arrest. \emph{New England Journal of Medicine} 351(7):637--646,
  \urlprefix\url{http://dx.doi.org/10.1056/NEJMoa040566}, pMID: 15306665.

\bibitem[{Hazinski et~al.(2005)Hazinski, Idris, Kerber, Epstein, Atkins, Tang,
  \protect\BIBand{} Lurie}]{Hazinski2005}
Hazinski MF, Idris AH, Kerber RE, Epstein A, Atkins D, Tang W, Lurie K (2005)
  Lay rescuer automated external defibrillator ({\textquotedblleft}public
  access defibrillation{\textquotedblright}) programs. \emph{Circulation}
  111(24):3336--3340, ISSN 0009-7322,
  \urlprefix\url{http://dx.doi.org/10.1161/CIRCULATIONAHA.105.165674}.

\bibitem[{Jones \protect\BIBand{} Henderson(2009)}]{Jones2009}
Jones M, Henderson D (2009) Maximum likelihood kernel density estimation: On
  the potential of convolution sieves. \emph{Computational Statistics and Data
  Analysis} 53(10):3726 -- 3733, ISSN 0167-9473,
  \urlprefix\url{http://dx.doi.org/https://doi.org/10.1016/j.csda.2009.03.019}.

\bibitem[{Kim et~al.(2017)Kim, Lim, Cho, \protect\BIBand{}
  C{\^o}t{\'e}}]{Kim2017}
Kim SJ, Lim GJ, Cho J, C{\^o}t{\'e} MJ (2017) Drone-aided healthcare services
  for patients with chronic diseases in rural areas. \emph{Journal of
  Intelligent {\&} Robotic Systems} 88(1):163--180, ISSN 1573-0409,
  \urlprefix\url{http://dx.doi.org/10.1007/s10846-017-0548-z}.

\bibitem[{Knight et~al.(2012)Knight, Harper, \protect\BIBand{}
  Smith}]{Knight2012}
Knight V, Harper P, Smith L (2012) Ambulance allocation for maximal survival
  with heterogeneous outcome measures. \emph{Omega} 40(6):918 -- 926, ISSN
  0305-0483,
  \urlprefix\url{http://dx.doi.org/https://doi.org/10.1016/j.omega.2012.02.003},
  special Issue on Forecasting in Management Science.

\bibitem[{Krishnan et~al.(2016)Krishnan, Marla, \protect\BIBand{}
  Yue}]{Krishnan2016}
Krishnan K, Marla L, Yue Y (2016) Robust ambulance allocation using risk-based
  metrics. \emph{2016 8th International Conference on Communication Systems and
  Networks (COMSNETS)}, 1--6, ISSN 2155-2509,
  \urlprefix\url{http://dx.doi.org/10.1109/COMSNETS.2016.7439958}.

\bibitem[{Kumar \protect\BIBand{} Michael(2012)}]{Kumar2012}
Kumar V, Michael N (2012) Opportunities and challenges with autonomous micro
  aerial vehicles. \emph{The International Journal of Robotics Research}
  31(11):1279--1291,
  \urlprefix\url{http://dx.doi.org/10.1177/0278364912455954}.

\bibitem[{Lerner et~al.(2012)Lerner, Rea, Bobrow, Acker, Berg, Brooks, Cone,
  Gay, Gent, Mears, Nadkarni, O{\textquoteright}Connor, Potts, Sayre, Swor,
  \protect\BIBand{} Travers}]{Lerner2012}
Lerner EB, Rea TD, Bobrow BJ, Acker JE, Berg RA, Brooks SC, Cone DC, Gay M,
  Gent LM, Mears G, Nadkarni VM, O{\textquoteright}Connor RE, Potts J, Sayre
  MR, Swor RA, Travers AH (2012) Emergency medical service dispatch
  cardiopulmonary resuscitation prearrival instructions to improve survival
  from out-of-hospital cardiac arrest. \emph{Circulation} ISSN 0009-7322,
  \urlprefix\url{http://dx.doi.org/10.1161/CIR.0b013e31823ee5fc}.

\bibitem[{Li et~al.(2011)Li, Zhao, Zhu, \protect\BIBand{} Wyatt}]{Li2011}
Li X, Zhao Z, Zhu X, Wyatt T (2011) Covering models and optimization techniques
  for emergency response facility location and planning: a review.
  \emph{Mathematical Methods of Operations Research} 74(3):281--310, ISSN
  1432-5217, \urlprefix\url{http://dx.doi.org/10.1007/s00186-011-0363-4}.

\bibitem[{Lin et~al.(2011)Lin, Morrison, \protect\BIBand{} Brooks}]{Lin2011}
Lin S, Morrison LJ, Brooks SC (2011) Development of a data dictionary for the
  strategies for post arrest resuscitation care (sparc) network for post
  cardiac arrest research. \emph{Resuscitation} 82(4):419 -- 422, ISSN
  0300-9572,
  \urlprefix\url{http://dx.doi.org/https://doi.org/10.1016/j.resuscitation.2010.12.006}.

\bibitem[{Marianov \protect\BIBand{} ReVelle(1994)}]{Marianov1994}
Marianov V, ReVelle C (1994) The queuing probabilistic location set covering
  problem and some extensions. \emph{Socio-Economic Planning Sciences}
  28(3):167 -- 178, ISSN 0038-0121,
  \urlprefix\url{http://dx.doi.org/https://doi.org/10.1016/0038-0121(94)90003-5}.

\bibitem[{Marianov \protect\BIBand{} ReVelle(1996)}]{Marianov1996}
Marianov V, ReVelle C (1996) The queueing maximal availability location
  problem: A model for the siting of emergency vehicles. \emph{European Journal
  of Operational Research} 93(1):110 -- 120, ISSN 0377-2217,
  \urlprefix\url{http://dx.doi.org/https://doi.org/10.1016/0377-2217(95)00182-4}.

\bibitem[{Marianov \protect\BIBand{} Serra(1998)}]{Marianov1998}
Marianov V, Serra D (1998) Probabilistic, maximal covering location allocation
  models for congested systems. \emph{Journal of Regional Science}
  38(3):401--424.

\bibitem[{Marianov \protect\BIBand{} Serra(2002)}]{Marianov2002}
Marianov V, Serra D (2002) Location--allocation of multiple-server service
  centers with constrained queues or waiting times. \emph{Annals of Operations
  Research} 111(1):35--50, ISSN 1572-9338,
  \urlprefix\url{http://dx.doi.org/10.1023/A:1020989316737}.

\bibitem[{Momont(2014)}]{Delft}
Momont A (2014) Ambulance drone. Technical report, Delft University of
  Technology,
  \urlprefix\url{https://www.tudelft.nl/en/ide/research/research-labs/applied-labs/ambulance-drone/}.

\bibitem[{Morrison et~al.(2008)Morrison, Nichol, Rea, Christenson, Callaway,
  Stephens, Pirrallo, Atkins, Davis, Idris, \protect\BIBand{}
  Newgard}]{Morrison2008}
Morrison L, Nichol G, Rea T, Christenson J, Callaway C, Stephens S, Pirrallo R,
  Atkins D, Davis D, Idris A, Newgard C (2008) Rationale, development and
  implementation of the resuscitation outcomes consortium epistry-cardiac
  arrest. \emph{Resuscitation} 78(2):161--169, ISSN 0300-9572,
  \urlprefix\url{http://dx.doi.org/10.1016/j.resuscitation.2008.02.020}.

\bibitem[{Mourelo~Ferrandez et~al.(2016)Mourelo~Ferrandez, Harbison, Weber,
  Sturges, \protect\BIBand{} Rich}]{Ferrandez2016}
Mourelo~Ferrandez S, Harbison T, Weber T, Sturges R, Rich R (2016) Optimization
  of a truck-drone in tandem delivery network using k-means and genetic
  algorithm. \emph{Journal of Industrial Engineering and Management} 9:374.

\bibitem[{Murray \protect\BIBand{} Chu(2015)}]{Murray2015}
Murray CC, Chu AG (2015) The flying sidekick traveling salesman problem:
  Optimization of drone-assisted parcel delivery. \emph{Transportation Research
  Part C: Emerging Technologies} 54:86 -- 109, ISSN 0968-090X,
  \urlprefix\url{http://dx.doi.org/https://doi.org/10.1016/j.trc.2015.03.005}.

\bibitem[{Pulver et~al.(2016)Pulver, Wei, \protect\BIBand{} Mann}]{Pulver2016}
Pulver A, Wei R, Mann C (2016) Locating aed enabled medical drones to enhance
  cardiac arrest response times. \emph{Prehospital Emergency Care}
  20(3):378--389,
  \urlprefix\url{http://dx.doi.org/10.3109/10903127.2015.1115932}, pMID:
  26852822.

\bibitem[{Reece(2014)}]{Reece}
Reece P (2014) Drones that save lives. Technical report, Drone Delivery Canada,
  \urlprefix\url{https://saltspringexchange.com/2014/06/26/drones-that-save-lives/}.

\bibitem[{Restrepo et~al.(2008)Restrepo, Henderson, \protect\BIBand{}
  Topaloglu}]{Restrepo2008}
Restrepo M, Henderson SG, Topaloglu H (2008) Erlang loss models for the static
  deployment of ambulances. \emph{Health Care Management Science} 12(1):67,
  ISSN 1572-9389, \urlprefix\url{http://dx.doi.org/10.1007/s10729-008-9077-4}.

\bibitem[{ReVelle \protect\BIBand{} Hogan(1988)}]{ReVelle1988}
ReVelle C, Hogan K (1988) A reliability-constrained siting model with local
  estimates of busy fractions. \emph{Environment and Planning B: Planning and
  Design} 15(2):143--152, \urlprefix\url{http://dx.doi.org/10.1068/b150143}.

\bibitem[{ReVelle \protect\BIBand{} Hogan(1989{\natexlab{a}})}]{ReVelle1989a}
ReVelle C, Hogan K (1989{\natexlab{a}}) The maximum availability location
  problem. \emph{Transportation Science} 23(3):192--200,
  \urlprefix\url{http://dx.doi.org/10.1287/trsc.23.3.192}.

\bibitem[{ReVelle \protect\BIBand{} Hogan(1989{\natexlab{b}})}]{ReVelle1989b}
ReVelle C, Hogan K (1989{\natexlab{b}}) The maximum reliability location
  problem and $\alpha$-reliablep-center problem: Derivatives of the
  probabilistic location set covering problem. \emph{Annals of Operations
  Research} 18(1):155--173, ISSN 1572-9338,
  \urlprefix\url{http://dx.doi.org/10.1007/BF02097801}.

\bibitem[{Rockafellar \protect\BIBand{} Uryasev(2000)}]{Rockafellar2000}
Rockafellar R, Uryasev S (2000) Optimization of conditional value-at-risk.
  \emph{Journal of Risk} 2(3):21--41.

\bibitem[{Rockafellar \protect\BIBand{} Uryasev(2002)}]{Rockafellar2002}
Rockafellar R, Uryasev S (2002) Conditional value-at-risk for general loss
  distributions. \emph{Journal of Banking and Finance} 26(7):1443--1471.

\bibitem[{Sch{\"o}llig et~al.(2011)Sch{\"o}llig, Hehn, Lupashin,
  \protect\BIBand{} D'Andrea}]{Schollig2011}
Sch{\"o}llig A, Hehn M, Lupashin S, D'Andrea R (2011) Feasiblity of motion
  primitives for choreographed quadrocopter flight. \emph{American Control
  Conference (ACC), 2011}, 3843--3849 (IEEE).

\bibitem[{Scott \protect\BIBand{} Scott(2017)}]{Scott2017}
Scott JE, Scott CH (2017) Drone delivery models for healthcare. \emph{HICSS}.

\bibitem[{Serra \protect\BIBand{} Marianov(1998)}]{Serra1998}
Serra D, Marianov V (1998) The p-median problem in a changing network: the case
  of barcelona. \emph{Location Science} 6:383--394.

\bibitem[{Sheather \protect\BIBand{} Jones(1991)}]{Sheather1991}
Sheather SJ, Jones MC (1991) A reliable data-based bandwidth selection method
  for kernel density estimation. \emph{Journal of the Royal Statistical
  Society. Series B (Methodological)} 53(3):683--690, ISSN 00359246,
  \urlprefix\url{http://www.jstor.org/stable/2345597}.

\bibitem[{Siddiq et~al.(2013)Siddiq, Brooks, \protect\BIBand{}
  Chan}]{Siddiq2013}
Siddiq AA, Brooks SC, Chan TC (2013) Modeling the impact of public access
  defibrillator range on public location cardiac arrest coverage.
  \emph{Resuscitation} 84(7):904 -- 909, ISSN 0300-9572,
  \urlprefix\url{http://dx.doi.org/https://doi.org/10.1016/j.resuscitation.2012.11.019}.

\bibitem[{Sun et~al.(2016)Sun, Demirtas, Brooks, Morrison, \protect\BIBand{}
  Chan}]{Sun2016}
Sun CL, Demirtas D, Brooks SC, Morrison LJ, Chan TC (2016) Overcoming spatial
  and temporal barriers to public access defibrillators via optimization.
  \emph{Journal of the American College of Cardiology} 68(8):836 -- 845, ISSN
  0735-1097,
  \urlprefix\url{http://dx.doi.org/https://doi.org/10.1016/j.jacc.2016.03.609}.

\bibitem[{Toro-D{\'\i}az et~al.(2013)Toro-D{\'\i}az, Mayorga, Chanta,
  \protect\BIBand{} McLay}]{ToroDiaz2013}
Toro-D{\'\i}az H, Mayorga ME, Chanta S, McLay LA (2013) Joint location and
  dispatching decisions for emergency medical services. \emph{Computers and
  Industrial Engineering} 64(4):917 -- 928, ISSN 0360-8352,
  \urlprefix\url{http://dx.doi.org/https://doi.org/10.1016/j.cie.2013.01.002}.

\bibitem[{Valenzuela et~al.(1997)Valenzuela, Roe, Cretin, Spaite,
  \protect\BIBand{} Larsen}]{Valenzuela1997}
Valenzuela TD, Roe DJ, Cretin S, Spaite DW, Larsen MP (1997) Estimating
  effectiveness of cardiac arrest interventions. \emph{Circulation}
  96(10):3308--3313, ISSN 0009-7322,
  \urlprefix\url{http://dx.doi.org/10.1161/01.CIR.96.10.3308}.

\bibitem[{Weisfeldt et~al.(2010)Weisfeldt, Sitlani, Ornato, Rea, Aufderheide,
  Davis, Dreyer, Hess, Jui, Maloney, Sopko, Powell, Nichol, \protect\BIBand{}
  Morrison}]{Weisfeldt2010}
Weisfeldt ML, Sitlani CM, Ornato JP, Rea T, Aufderheide TP, Davis D, Dreyer J,
  Hess EP, Jui J, Maloney J, Sopko G, Powell J, Nichol G, Morrison LJ (2010)
  Survival after application of automatic external defibrillators before
  arrival of the emergency medical system: Evaluation in the resuscitation
  outcomes consortium population of 21 million. \emph{Journal of the American
  College of Cardiology} 55(16):1713 -- 1720, ISSN 0735-1097,
  \urlprefix\url{http://dx.doi.org/https://doi.org/10.1016/j.jacc.2009.11.077}.

\bibitem[{Xia et~al.(2017)Xia, Batta, \protect\BIBand{} Nagi}]{Xia2017}
Xia Y, Batta R, Nagi R (2017) Controlling a fleet of unmanned aerial vehicles
  to collect uncertain information in a threat environment. \emph{Operations
  Research} 65(3):674--692,
  \urlprefix\url{http://dx.doi.org/10.1287/opre.2017.1590}.

\bibitem[{Zayas-Cab{\'a}n et~al.(2013)Zayas-Cab{\'a}n, Lewis, Olson,
  \protect\BIBand{} Schmitz}]{Zayas2013}
Zayas-Cab{\'a}n G, Lewis ME, Olson M, Schmitz S (2013) Emergency medical
  service allocation in response to large-scale events. \emph{IIE Transactions
  on Healthcare Systems Engineering} 3(1):57--68,
  \urlprefix\url{http://dx.doi.org/10.1080/19488300.2012.762816}.

\end{thebibliography}


\end{document}